\documentclass[12pt]{article}
\usepackage{amsmath,amssymb,amsthm,setspace}
       \newtheorem{theorem}{Theorem}
       
       \newtheorem{lemma}[theorem]{Lemma}

       \theoremstyle{definition}
       \newtheorem{definition}{Definition}
       \newtheorem{remark}{Remark}
       
\def\N{\ensuremath{\mathbb N}} 
\def\R{\ensuremath{\mathbb R}} 
\def\Z{\ensuremath{\mathbb Z}} 
\newcommand{\aaa}{\alpha}

\newcommand{\cad}{{\cal D}}

\newcommand{\caj}{{\cal J}}
\newcommand{\cak}{{\cal K}}
\newcommand{\call}{{\cal L}}
\newcommand{\capp}{{\cal P}}

\newcommand{\cas}{{\cal S}}

\newcommand{\cat}{{\cal T}}

\newcommand{\ddd}{\delta}

\newcommand{\defeq}{{\buildrel {\rm def}\over =}}
\newcommand{\ds}{\displaystyle}

\newcommand{\eqq}{\equiv}
\newcommand{\ess}{\emptyset}
\newcommand{\FFF}{\Phi}

\newcommand{\gggg}{\gamma}
\newcommand{\GGG}{\Gamma}

\newcommand{\hFFF}{\widehat{\Phi}}

\newcommand{\hPPP}{\widehat{\Psi}}

\newcommand{\kkk}{\kappa}
\newcommand{\lf}{\lfloor}
\newcommand{\LLLg}{\Lambda_{\gamma }}
\newcommand{\llll}{\lambda}
\newcommand{\LLL}{\Lambda}

\newcommand{\oaaa}{{\overline{\alpha}}}
\newcommand{\ocap}{{\overline{{\cal P}}}}
\newcommand{\oEE}{{\overline{E}}}
\newcommand{\off}{{\overline{f}}}

\newcommand{\ointt}{{\overline{\int}}}
\newcommand{\oLLLg}{{\overline{\Lambda}}_{\gamma }}
\newcommand{\oLLL}{{\overline{\Lambda}}}
\newcommand{\olll}{{\overline{\lambda}}}
\newcommand{\oooo}{{\overline{\omega}}}
\newcommand{\OOO}{\Omega}
\newcommand{\oo}{\infty}
\newcommand{\ottt}{{\overline{\tau}}}
\newcommand{\oXXX}{\overline{\Xi}}
\newcommand{\oXX}{{\overline{X}}}
\newcommand{\oyy}{{\overline{y}}}

\newcommand{\PPP}{\Psi}

\newcommand{\rf}{\rfloor}
\newcommand{\sm}{\setminus}
\newcommand{\sse}{\subset}
\newcommand{\tCC}{\widetilde{C}}
\newcommand{\tEE}{{\widetilde{E}}}
\newcommand{\tFFF}{{\widetilde{\Phi}}}
\newcommand{\tPPP}{{\widetilde{\Psi}}}
\newcommand{\tqq}{\widetilde{q}}
\newcommand{\trrr}{{\widetilde{\rho}}}
\newcommand{\TTT}{\Theta}
\newcommand{\tttt}{\tau}

\newcommand{\tUU}{\widetilde{U}}

\newcommand{\veee}{{\varepsilon}}

\title{Divergent Square Averages}
\author{
Zolt\'an Buczolich\thanks{This
author's work was started during his
 visit at the Department of Mathematics of University of
North Texas and later was supported by the Hungarian
National Foundation for Scientific Research T049727.},
Department of Analysis, E\"otv\"os Lor\'and\\
University, P\'azm\'any P\'eter S\'et\'any 1/c, 1117 Budapest, Hungary\\
email: buczo@cs.elte.hu\\
{\tt www.cs.elte.hu/\hbox{$\sim$}buczo}\\
 \\ and\\
 \\
R. Daniel Mauldin\thanks{Supported in part by NSF grants DMS 0400481, 0700831
and as a van Fleck Visiting Researcher at Wesleyan
University during 2005.
}, Department of Mathematics,\\
 University of North Texas, Denton, Texas 76203-1430, USA\\
email: mauldin@unt.edu\\
{\tt www.math.unt.edu/\hbox{$\sim$}mauldin} }
\date{\today}
\begin{document}
\maketitle
\begin{abstract}
In this paper we answer a question of  
J. Bourgain which was motivated by questions
A. Bellow and H. Furstenberg. We show
that the sequence $\{ n^{2}\}_{n=1}^{\oo}$
 is $L^{1}$-universally bad.
This implies that it is not true that
given a dynamical system $(X ,\Sigma, \mu, T)$ and $f\in L^{1}(\mu)$,
the ergodic means
$$
\lim_{N\to \oo}\frac{1}N\sum _{n=1}^{N}f(T^{n^{2}}(x))
$$
converge
almost surely.
\end{abstract}

\section{Introduction}
Research related to almost everywhere convergence of 
ergodic averages along the squares
was initiated by questions of
 A. Bellow (see \cite{OW}) and of H. Furstenberg \cite{FU}.
Results of Bourgain \cite{B1}, \cite{B2}, \cite{B3} imply that
given a dynamical system $(X ,\Sigma, \mu, T)$ and $f\in L^{p}(\mu)$,
for some $p>1$, the ergodic means
\begin{equation}\label{bourg1}
\lim_{N\to \oo}\frac{1}N\sum _{n=1}^{N}f(T^{n^{2}}(x))
\end{equation}
converge
almost surely.
Bourgain also asked in \cite{B1}, \cite{B5} whether this result is true
for $L^{1}$ functions.
In this paper we give a negative answer to this question.

Let us recall some concepts related to this problem.
\begin{definition}
A sequence $\{ n_{k}\}_{k=1}^{\oo}$ is $L^{1}$-{\it universally bad}
if for all ergodic dynamical systems there is some
$f\in L^{1}$ such that
$$\lim_{N\to\oo}\frac{1}{N}\sum_{k=1}^{N}
f(T^{n_{k}}x)$$
fails to exist for all $x$ in a set of positive measure.
\end{definition}

By the Conze principle and the Banach principle
of Sawyer (see \cite{[con]},
\cite{[saw]}, or \cite{[RW]}),
the sequence $\{ n_{k}\}_{k=1}^{\oo}$ is $L^{1}$-universally
bad if there is no constant $C<\oo$ such that for all
systems $(X,\Sigma,\mu,T)$
 and all $f\in L^{1}(\mu)$
we have the following weak $(1,1)$ inequality
for all $t\in \R:$
\begin{equation}\label{weaklolo}
\mu\left ( \left \{ x:\sup_{N\geq 1}\left |
\frac{1}N
\sum_{k=1}^{N}f(T^{n_{k}}x)\right|>t
\right \}\right)\leq \frac{C}{t}\int |f| d\mu.
\end{equation}

The main result of this paper is
\begin{theorem}\label{mainth}
The sequence $\{ k^{2}\}_{k=1}^{\oo}$ is $L^{1}$-universally bad.
\end{theorem}

This theorem will be proved by showing that there is no
constant $C$ such that the weak $(1,1)$ inequality given in
\eqref{weaklolo} holds.

This paper is a new and substantially modified
version of our preprint from 2003. The proof
in that preprint contained a gap but the methods
of that paper lead to a solution
of a counting problem raised by
I. Assani (see \cite{[abm]} and \cite{[abm2]}).

The paper is organized as follows.
In Section \ref{nth} we develop the necessary ingredients
we need concerning the asymptotic distribution of squares
modulo $q$ where $q$ is the product of $\kappa$ distinct
primes. The specific technical property we need is
given in Lemma \ref{probabilitylemma}. In Section
\ref{perr} we develop the notion of a periodic
rearrangement of a given periodic set. Lemma
\ref{taurearrangement} states a property about the
frequency squares hit such sets;
we will need this later in our construction.
Section \ref{secksystem} is the technical heart of
the paper. For positive integers $K$, $M$ and a periodic set
$\Lambda$ we define the notion of a $K-M$ family living
on $\Lambda$. What we need for the proof of our main theorem
is the existence of some specific families living
on $\Lambda=\R.$ The properties of these families are stated
in Lemma \ref{lksystem}. However, we need a double induction
argument to show that such families exist.
In Section \ref{41}, in Lemma \ref{lputres} assuming $K-M$ families
exist for all parameter values on $\R$, we show that they exist
on periodic sets $\LLL$. In Section \ref{42} we turn to the
proof that if $K-M$ families exist, then $(K+1)-M$ families
exist as well, this induction on $K$ is our outer inductive construction.
In Sections \ref{421} through \ref{428} we carry out the first
step of this induction, while in Sections \ref{422b} through \ref{428b}
 we show
how this first step of the induction should be altered for $(K+1)-M$
families when $K>0$.
The proof of the existence of $(K+1)-M$ families involves an intricate
inner inductive construction, the ``leakage process", which
is outlined in Section \ref{421} and carried out in Sections
\ref{422} through \ref{427}. Once it has halted, it is
shown in Section \ref{428} how to adjust the functions
so that the next stage of the outer induction holds.
In Section \ref{secmain}, we give the proof of the
main theorem. We construct a sequence of rational
rotations $T_{p}$, functions $f_{p}$ and numbers $t_{p}$
which witness that there is no
constant $C$ satisfying \eqref{weaklolo}.

To understand the heuristics behind our proof
it might also be useful to look at \cite{BM}.

 If someone prefers to have a general overview of the
main ideas of the paper before turning to the details here
is a recommended quick tour: After reading the introduction
read Definition \ref{conmkilenc}.
Then jump to Section \ref{secmain} and read the statement of
Theorem \ref{*c1c}. Skip the proof of this theorem
and read the details of the proof of Theorem
\ref{mainth} at the end of Section \ref{secmain}.
Then continue with Section \ref{nth} and read it until
\eqref{*3bud}. Jump to Definition \ref{5dfw} and
Remark \ref{**5dfw}. Then read Section \ref{secksystem}
starting at Definition \ref{nagykompl}
until the statements of Lemmas \ref{lksystem} and
\ref{lputres}. From here jump to Subsection \ref{42}
and read it until the paragraph
above
\ref{421}. Continue with Lemma \ref{probabilitylemma}
and the paragraph above it. Read Remark \ref{55dfw}. Then jump
to Section \ref{perr} and read it until Lemma \ref{taurearrangement}
is stated. Finally, read Subsection \ref{421} until the paragraph
containing \eqref{*h6}.

Let us fix some notation.
Given $f:\R\to\R$, periodic by $p$ we put
$$\ointt f=\frac{1}p \int_{0}^{p}f(x)dx.$$
Given a Lebesgue measurable set $A$, periodic by
$p$ we put $$\olll(A)=\frac{1}p\lambda (A\cap [0,p))=\lim_{N\to\oo}\frac{\llll(A\cap [-N,N])}
{2N}.$$

\section{Number Theory/Quadratic residues}\label{nth}

For each $q\in\N$ and $n\in \Z$ set $\veee(n,q)=1$
if $n$ is congruent to a square modulo $q$,
and let $\veee(n,q)=0$ if not.
We denote by $\sigma _{q}$ the number of
squares modulo $q$.
If $p$ is an odd prime, then $\sigma _{p}=\frac{p+1}2$.
If $q=p_{1}\cdot \cdot \cdot p_{\kappa }$ where
$p_{1},...,p_{\kappa }$ are distinct odd primes,
then (by the fact that something is a square modulo
$q$ if and only if it is a square modulo
each $p_{i}$ plus by using the Chinese remainder
theorem)
$\sigma _{q}=\prod_{i=1}^{\kappa }\frac{p_{i}+1}2$.
For elementary properties of quadratic residues
see \cite{HW} pages 67-69, or Chapter 3 of
\cite{NZ}. We remark that though
$0^2=0$, when talking about quadratic residues
usually only those are considered which are not
congruent to $0$, but since $\veee(n,q)=1$
when $n$ is congruent to $0$
modulo $q$
we will regard
$0$ a quadratic residue (or square) in this paper.

In Section \ref{perr} we will use the {\it Legendre symbol}.
If $\tttt$ is an odd prime the Legendre symbol $\ds 
 \bigg(\frac{n}{\tttt}\bigg)=
\bigg(\frac{n}{\tttt}\bigg)_{\call}$
equals $0$ if $\tttt$ divides $n$, otherwise it equals
$+1$ if $n$ is a square modulo
$\tttt$ and equals $-1$ if $n$ is not a square modulo $\tttt$.
To avoid notational confusion we will use the subscript $\call$
for the Legendre symbol.
It is a character, that is, $\ds \bigg(\frac{n m}{\tttt}\bigg)_{\call}=
\bigg(\frac{n}{\tttt}\bigg)_{\call}\bigg(\frac{m}{\tttt}\bigg)_{\call}$,
and $\ds 1+\bigg(\frac{n}{\tttt}\bigg)_{\call}$ equals the number of solutions
$x$ $mod$ $\tttt$ to the congruence $x^{2}\eqq n$ $mod$ $\tttt$.

Put $\LLL_{0}(q)= \{ n\in \Z: \veee(n,q)=1  \}$, the set of 
integers which are quadratic residues
modulo $q$.
Clearly,
\begin{equation}\label{*3bud}
\#(\LLL_{0}(q)\cap [0,q))=\sigma_{q}>\frac{q}{2^{\kappa }}.
\end{equation}

Next we discuss some results from \cite{[Per]} concerning
the distribution of the squares modulo $q$.
Given $K$, consider a fixed sequence $(\veee_{1},...,\veee_{K})$
of zeros and ones. Assume that $a_{1},...,a_{K}$ are distinct
integers modulo an odd prime $p$.

Set $\nu_{p}(\veee_{1},...,\veee_{K};a_{1},...,a_{K})=
\#  \{ n:0\leq n<p,\, \veee(n+a_{i},p)=\veee_{i}
\text{ for } i=1,...,K  \},$
that is, $\nu_{p}(\veee_{1},...,\veee_{K};
a_{1},...,a_{K})$ counts the number of occurrences
of $(\veee_{1},...,\veee_{K})$ in translated copies of
$a_{1},...,a_{K}$ modulo $p$. Then
\begin{equation}\nonumber
\frac{p}{2^{K}}-K(3+\sqrt p)\leq
\nu_{p}(\veee_{1},...,\veee_{K};a_{1},...,a_{K})\leq
\frac{p}{2^{K}}+K(3+\sqrt p).
\end{equation}

The ``probability" of the occurrence of $(\veee_{1},...,\veee_{K})$
in translated copies of $(a_{1},...,a_{K})$
is
$$ \capp_{p}
(\veee_{1},...,\veee_{K};a_{1},...,a_{K})=\frac{\nu_{p}
(\veee_{1},...,\veee_{K};a_{1},...,a_{K})}p$$
and by the above result if $(a_{1},...,a_{K})$ is fixed,
then
\begin{equation}\label{*a2}
\capp_{p}(\veee_{1},...,\veee_{K};a_{1},...,a_{K})
\to \frac{1}{2^{K}}\text{ as the odd prime }p\to\oo.
\end{equation}

Next we want to choose square free numbers
$q=p_{1}\cdot \cdot \cdot p_{\kappa }$, where
$p_{1},...,p_{\kappa }$ are distinct sufficiently
large odd primes with good statistical properties.

A number $n$ is a square modulo
$q$ if and only if it is a square modulo
each of the primes $p_{1},...,p_{\kappa }.$
By \eqref{*a2}
given $a_{1},...,a_{K}$ and
keeping $\kappa $, fixed
as $\min \{ p_{1},...,p_{\kappa }  \}\to\oo$
we have
\begin{equation}\label{*a3}
\capp_{q}(\veee_{1},...,\veee_{K};a_{1},...,a_{K})
\to \left (\frac{1}{2^{\kappa }}\right )^{\sum_{i=1}^{K}\veee_{i}}
\left (1-\frac{1}{2^{\kappa }}\right)^{K-\sum_{i=1}^{K}\veee_{i}},
\end{equation}
that is, statistically squares modulo
$q$ look like outcomes of
independent Bernoulli trials
with probabilities $\frac{1}{2^{\kappa }}$ and $\left
(1- \frac{1}{2^{\kappa }}
\right).$
Without going into technical details of this fact
from number theory, we just give an outline of a proof by induction on
$\kappa $. For $\kappa =1$, \eqref{*a3}
follows from \eqref{*a2}. Suppose $\kappa > 1$ and \eqref{*a3} holds for
$\kappa -1$.  Set $q_{0}=p_{1}\cdot \cdot \cdot p_{\kappa -1}.$
For any possible
choice of $\veee' = (\veee_{1}',...,\veee_{K}')$ and $\veee'' = (\veee_{1}'',...,\veee_{K}'')$,
set
$W(\veee',\veee'')= \{ n:0\leq n<q ,\  \veee(n+a_{i},q_{0})=\veee_{i}'\text{ and
 }\veee(n+a_{i},p_{\kappa })=\veee_{i}''
\text{ for }i=1,...,K
\}$,
$W_{0}(\veee') = \{ n:0\leq n<q_0 ,\  \veee(n+a_{i},q_{0})=\veee_{i}'
\text{ for }i=1,...,K
\},$ and
$W_{\kappa }(\veee'')= \{ n:0\leq n<p_{\kappa } ,\
\veee(n+a_{i},p_{\kappa })=\veee_{i}''
\text{ for }i=1,...,K
\}.$
Observe that for any $n$ the numbers $n+jq_{0}$,
$j=0,...,p_{\kappa }-1$ hit
each residue class modulo $p_{\kappa }$ exactly once,
and $\veee(n+jq_{0}+a_{i},q_{0})=\veee(n+a_{i},q_{0})$
for all $i$.
Using this one can see that $W_{0}(\veee')$ and $W_{\kappa }(\veee'')$
are independent in the sense that $\#W(\veee',\veee'')=\#W_{0}(\veee')\#W_{\kappa }(\veee'')$.
For $\veee = (\veee_1,...,\veee_K)$, set $G(\veee) =  \{ n:0\leq n<q ,\
\veee(n+a_{i},q)=\veee_{i}
\text{ for }i=1,...,K
\}.$
Then $\# G(\veee) =  \sum \# W(\veee',\veee'')$, where the sum is taken over all pairs $(\veee',\veee'')$ whose
coordinatewise product is $\veee$. Taking the limit as $\min{p_1,...,p_\kappa}$ goes to infinity of $\# G(\veee)/q$,
using \eqref{*a2} for the limit of $\# W_{\kappa}(\veee'')/p_{\kappa}$ and \eqref{*a3} for the limit of $\# W_0(\veee')/q_0$,
and noting that the only
thing that matters in these limiting probabilites is the number of
$1$'s in the sequence, we have after setting $m = \sum_{i=1}^K \veee_i$,
that the limiting value is
$$
{\frac{1}{2^K}}\sum_{i =0}^{K-m}\binom{K-m}{i}
(\frac{1}{2^{\kappa-1}})^m (1-\frac{1}{2^{\kappa - 1}})^i =  (\frac{1}{2^{\kappa }})^{m}
(1-\frac{1}{2^{\kappa }})^{K-m},
$$
which is what we wanted.

Consider
an infinite sequence of pairwise independent
random variables
$X_{i}:Y\to  \{ 0,1  \}$ with $P(X_{i}(\omega )=1)=\frac{1}{2^{\kappa }}$,
$P(X_{i}(\omega )=0)=1-\frac{1}{2^{\kappa }}.$
Then $E(X_{i}(\omega ))=\frac{1}{2^{\kappa }}.$ By the law of large
numbers if $\rho >0$ and $K \to\oo$, then
$$P\left (\left | \frac{1}K\sum_{i=1}^{K}
X_{i}(\omega )-\frac{1}{2^{\kappa }}\right|\geq\rho\right)\to 0.$$
Given $\rho ,\epsilon_1 >0$ if $K$ is large enough, then
$$P\left (\left |\frac{1}K\sum_{i=1}^{K}
X_{i}(\omega )-\frac{1}{2^{\kappa
}}\right|\geq\rho\right)<\epsilon_1 .$$

For odd primes $p_{1}<p_{2}<...<p_{\kappa }$
put  $q=p_{1}\cdot \cdot \cdot p_{\kappa}.$
Given
distinct integers $a_{1},...,a_{K}$
consider
$X_{i}(n)=\veee(n+a_{i},q)$. By \eqref{*a3} as $p_{1}\to\oo$
the variables $X_{i}(n)$  approximate independent random
variables with Bernoulli distribution $\frac{1}{2^{\kappa
}}$, $\left (1-\frac{1}{2^{\kappa }}\right)$.
Hence,
given a sufficiently large $K$
if $p_{1}$ is sufficiently large then
$$\frac{\#  \{ n\in [0,q):\left  |
\frac{1}K \sum_{i=1}^{K}\veee(n+a_{i},q)-\frac{1}{2^{\kappa
}}\right| \geq \rho  \}}
{q}<\epsilon_1 .$$

In particular, we have
\begin{equation}\label{*4dfw}
\#  \{ n\in [0,q):
\sum_{i=1}^{K}\veee(n+a_{i},q)\geq \frac{K}{2^{\kappa }}+K\rho   \}<
\epsilon_1 q.
\end{equation}

For later arguments we now introduce some parameters.
We will need a suitably small
``leakage constant" $\gamma \in (0,1)$
of the form $\gamma =2^{-c_{\gamma }}$ where
$c_{\gamma }\in \N.$ Then we work with
$\kappa >c_{\gamma }.$
Furthermore we use small constants
$\rho,\rho _{1} >0.$
For large $K_{1}$ we have
\begin{equation}\label{**3}
K_{2}\defeq\lceil (1+\rho_{1} 2^{\kappa })\gamma K_{1}\rceil
<(1+\rho _{1}) (1+\rho_{1} 2^{\kappa })\gamma K_{1}
\text{ and }\gamma 2^{\kappa }/K_{1}<\epsilon _{1}.
\end{equation}
Thus for a large $K_{1}$ we can
choose $p_{1}'$ such that
\begin{equation}\label{*5budd}
p_{1}'>\max \{ K_{1}+\gamma 2^{\kappa }, K_{1}/\epsilon_1  \}
\end{equation}
and, in addition
if $p_{1}>p_{1}'$ we have for any $q=p_{1}\cdot \cdot \cdot p_{\kappa },$
$p_{1}<...<p_{\kappa }$,
\begin{equation}\label{*a7}
\#  \{ n\in [0,q):
\left |\frac{1}{K_{1}}
\sum_{i=1}^{K_{1}}\veee(n+i,q)-\frac{1}{2^{\kappa }}\right|\geq\rho   \}
<\epsilon_{1} q.
\end{equation}
Also, we can choose $p_{1}''\geq p_{1}'$
such that if $p_{1}>p_{1}''$ then
given any integers $a_{1},...,a_{K_{2}}$
so that the difference of any two of them is less than
$p_{1}'$ we also have, after a simple change of notation
in \eqref{*4dfw}
\begin{equation}\label{**a7}
\#  \{ n\in [0,q):
\sum_{i=1}^{K_{2}}\veee(-(n+a_{i}),q)\geq
 \frac{K_{2}}{2^{\kappa }}+K_{2}\rho_{1}
   \}
<\epsilon_1 q,
\end{equation}
the negative sign in the first argument of $\veee(\cdot ,q)$
is due to technical reasons in later arguments;
it is clear that
if $n$ takes all possible values modulo
$q$ then so does $-n.$

Let $n_{1}\in [0,\oo)$ be the first number for which
$$\cad(K_{1},n_{1},q)\defeq
\left |\frac{1}{K_{1}}\sum_{i=1}^{K_{1}}
\veee(n_{1}+i,q)-\frac{1}{2^{\kappa }}
\right|< \rho .$$

Next choose the least $n_{2}\geq n_{1}+K_{1}$
such that $\cad(K_{1},n_{2},q)< \rho .$

Continue and
set $\caj= \{ j\in \N: 0\leq n_{j}<q  \}.$
If $n'\in [0,q)\sm\cup_{j\in \caj}[n_{j},n_{j}+K_{1})$,
then $\cad(K_{1},n',q)\geq \rho $ holds.

Hence, by \eqref{*a7}
\begin{equation}\label{*a7b}
\#  \{ n\in [0,q):
n\not\in \bigcup_{j\in\caj}[n_{j},n_{j}+K_{1})  \}<
\epsilon_1 q.
\end{equation}
If $j\in \caj$ then by the definition of $\veee(n,q)$
\begin{equation}\label{*a8}
K_{1}(\frac{1}{2^{\kappa }}-\rho )
<\#\bigg( \LLL_{0}(q)\cap [n_{j}+1,n_{j}+K_{1}]
\bigg)
<K_{1}(\frac{1}{2^{\kappa }}+\rho ),
\end{equation}
so that the number of quadratic residues
modulo $q$ in the interval $[n_{j}+1,n_{j}+K_{1}]$
is approximately $K_{1}/{2^{\kappa }}.$

\begin{definition}\label{5dfw}
Set $\LLLg(q)=-\LLL_{0}(q)+ \{ j\in\Z:
0\leq j<\gamma 2^{\kappa }  \},$
$\oLLLg(q)=\LLLg(q)+[0,1)=
-\LLL_{0}(q)+ \{ x:0\leq x< \gamma 2^{\kappa }  \}.$
For ease of notation in the sequel if we have a fixed $\gamma $
and we do not want to emphasize the dependence on $\gamma $
we will write $\LLL(q)$ and $\oLLL(q)$, instead of
$\LLLg(q)$ and $\oLLLg(q)$, respectively.
(To make it easier to memorize our notation for these
$\LLL$ type sets it will be useful to keep in mind that the sets
$\LLL$ without the bars will be subsets of $\Z$ and  the sets $\oLLL\sse \R$
will be obtained from the corresponding $\LLL$ sets by adding $[0,1)$.
\end{definition}

\begin{remark}\label{**5dfw}

Here are some ``heuristic" comments related to the above definition.

If $\LLL_{0}(q)$ equaled $ \{ k\cdot 2^{\kappa }:k\in \Z  \}$,
then $\olll(\oLLL(q))$ would equal $\gamma .$

Next suppose that $\LLL_0(q)$ is the set of quadratic
residues.
If the intervals making up  $\oLLL_\gamma(q))$ were
disjoint, then  $\olll(\oLLL_\gamma(q))$ would be $\gamma\prod_{i=1}^\kappa
\frac{p_i+1}{p_i}$, somewhat larger than $\gamma$. However, by results in
\cite{[Kur]}
for a fixed $\gamma $ as $\kappa $ goes to $\oo$,
the normalized
 gaps between consecutive elements
of $\LLL_{0}(q)$ converge to an exponential distribution. We will
make explicit use of this fact in Lemma \ref{rdmlem}. 
Since the normalizing factor is
$\sigma_{q}$,
the number
of squares modulo $q$ approximately
equals $q/2^{\kappa }$ and the average value
of the spacing between elements of $\LLL_{0}(q)$
is very close to $2^{\kappa }$. For each $\kappa $ sufficiently large
if $p_{1}$  is sufficiently large, for a certain percentage  of different elements $i,i'\in \LLL_{0}(q)$ the intervals
$[i,i+\gamma 2^{\kappa })$ and $[i',i'+\gamma 2^{\kappa })$
will overlap. Under these conditions $\olll(\oLLL(q))$ will be less than
$\gamma $, but the smaller $\gamma $ is, the closer
$\olll(\oLLL(q))/\gamma $ is to $1$ for large $q$'s. We will
take advantage of this property particularly when we fix a leakage
constant.
\end{remark}

By Remark \ref{**5dfw}, 
 $\#((\Z\cap [0,q))\sm \LLL(q))$
is a little larger than $(1-\gamma )\#(\Z\cap [0,q))$ for large $\kappa$ and $p_1$.

Suppose $\rho >0$ is given. If a $q$-periodic set $A\sse \Z$
is ``not sufficiently $\rho $-independent" from $\LLL(q)$
it may happen that
\begin{equation}\label{rhoind}
\#((A\sm\LLL(Q))\cap[0,q))<(1-\rho )(1-\gamma )\#(A\cap[0,q)).
\end{equation}
In the next lemma we consider translated copies of $\LLL_{0}(q)$.
We show that for $q$'s with large prime factors
 there is only a small portion of $n$'s when $A_{n}=n+\LLL_{0}(q)$
 satisfies  \eqref{rhoind}.

\begin{lemma}\label{probabilitylemma}
Given
a positive integer
$\kappa $
 and
$\epsilon $, $\rho >0$ there exists
$p_{1}''$ such that if
the odd primes satisfy
$p_{1}''<p_{1}<...<p_{\kappa }$
and $q=p_{1}\cdot \cdot \cdot p_{\kappa }$,
then
\begin{align}\label{*a9}
\# \bigg\{ n\in [0,q):&\#\bigg(\big((n+\LLL_{0}(q))\sm
\LLL(q)\big)\cap [0,q)\bigg)<\\
\nonumber
&(1-\rho )(1-\gamma )\#\bigg(\LLL_{0}(q)\cap
[0,q)\bigg) \bigg \}<\epsilon q.
\end{align}
\end{lemma}

\begin{remark}\label{55dfw}
The heuristics behind \eqref{*a9} are the following.
The number of $n\in [0,q)$ for which $n\in \LLL(q)$
is a little less than $\gamma q$,
due to ``overlaps".
This means that the number of those $n$'s for which
$n\not\in \LLL(q)$ is a little larger than
$(1-\gamma )q.$ Now one can examine what happens
when we look at translated copies of
$\LLL_{0}(q).$ Formula \eqref{*a9} says
that for ``most" translated copies of $\LLL_{0}(q)$
we cannot have much less than $(1-\gamma )\#
(\LLL_{0}(q)\cap [0,q))$ elements of $n+\LLL_{0}(q)$
outside $\LLL(q).$

Before beginning the proof of Lemma \ref{probabilitylemma}
we choose $\rho _{0}>0$ such that
\begin{equation}\label{*6budd}
\left ( 1-\frac{1}{(1-\rho _{0})^{2}}\gamma \right)
>(1-\rho )(1-\gamma ).
\end{equation}
Recall
from number theory
that if $p_{1}''$ is sufficiently large then we have
\begin{equation}\label{*5dfw}
(1-\rho_{0} )\frac{q}{2^{\kappa }}<\#(\LLL_{0}(q)
\cap [0,q))=\sigma _{q}
<\frac{1}{(1-\rho_0 )}\frac{q}{2^{\kappa }}.
\end{equation}
\end{remark}

\begin{proof}
We can assume that $\epsilon <1.$
Take $0<\epsilon _{1}< \epsilon ^{2}/32$
and then choose
$K_{1}$, $p_{1}'$ and $p_{1}''$
as above.
By \eqref{*5budd} and $q>p_{1}''\geq p_{1}'$ we have
$q>K_{1}/\epsilon _{1}$ and
by \eqref{*a7b}
\begin{equation}\label{*a10}
\# \{ n\in [0,q):
n\not \in\bigcup_{j\in\caj}[n_{j},n_{j}+K_{1})  \}
<{\epsilon_{1} }q
\end{equation}
for a suitable index set $\caj$,
defined above,
and for each $j\in\caj$
we have $\cad(K_{1},n_{j},q)< \rho _{1}$
for a suitable $\rho _{1}>0.$
This means that
\begin{equation}\label{**a10}
K_{1}\left (\frac{1}{2^{\kappa }}-\rho _{1}\right)<
\#(\LLL_{0}(q)\cap
[n_{j}+1,n_{j}+K_{1}])<
K_{1}\left (\frac{1}{2^{\kappa }}+\rho _{1}\right).
\end{equation}
By Definition
\ref{5dfw} for $n\in [0,q)$
we have $n+n'\in \LLL(q)=\LLL_{\gggg}(q)$
if $n+n'-i\in -\LLL_{0}(q)$
holds for an $i=0,...,\gamma 2^{\kappa }-1.$
Set
$\LLL_{0}^{j}(q)=\LLL_{0}(q)\cap[n_{j}+1,n_{j}+K_{1}]
$
and
$$\LLL^{j}(q)=-\left
(\bigcup_{i=0}^{\gamma 2^{\kappa }-1}(\LLL_{0}^{j}(q)-i)\right )$$
for each $j\in \caj.$ Using the definition
of $K_{2}$ in \eqref{**3},
and \eqref{**a10}
choose distinct numbers $a_{i',j}$, $i'=1,...,K_{2}$,
so that
\begin{equation}\label{*9zz}
-\LLL^{j}(q)=
\bigcup_{i=0}^{\gamma 2^{\kappa }-1}
\LLL_{0}^{j}(q)-i
\subseteq A_{j}\defeq
\bigcup_{i'=1}^{K_{2}} \{ a_{i',j}  \}
\subset
[n_{j}-\gamma 2^{\kappa },n_{j}+K_{1}]
.
\end{equation}
Clearly, by the choice of $p_{1}'$
in \eqref{*5budd}
the difference of any two of the
$a_{i',j}$'s
 is less
than $p_{1}'$.
By \eqref{*9zz}
$$(n+\LLL_{0}(q))\cap \LLL^{j}(q)
\sse (n+\LLL_{0}(q))\cap (-A_{j}).$$

Observe that there exists $n'\in \LLL_{0}(q)$
such that $n+n'\in -A_{j}$ if and only if
there exists $i'$ such that $n+n'=-a_{i',j}$,
that is, $n+a_{i',j}=-n'\in -\LLL_{0}(q).$

Recall that $n+a_{i',j}\in -\LLL_{0}(q)$ if and only if
$\veee(-(n+a_{i',j}),q)=1.$
Set $N_{j}= \{ n\in [0,q):
\sum_{i'=1}^{K_{2}}\veee(-(n+a_{i',j}),q)\geq
K_{2}(\frac{1}{2^{\kappa }}+\rho _{1})  \}.$
If $n\not \in N_{j}$, $n\in [0,q)$ then
\begin{equation}\label{*6x}
\#\bigg((n+\LLL_{0}(q))\cap \LLL^{j}(q)\bigg)
< K_{2}(\frac{1}{2^{\kappa }}
+\rho _{1}).
\end{equation}
By \eqref{**a7}
\begin{equation}\label{**5}\# N_{j}<\epsilon _{1}q.\end{equation}
Here we remark that \eqref{**a7} can be used so that
we have \eqref{**5} for all $j\in\caj.$
Indeed, assume $A\sse [-\gamma 2^{\kappa },K_{1}]$ and
$A+n_{j}=A_{j}$ then
$$N_{j}=N_{A}\defeq
 \{ n\in [0,q):
\sum_{a\in A}\veee(-(n+a),q)\geq K_{2}(\frac{1}{2^{\kappa }}+\rho _{1})  \}.$$
There are $2^{\gamma 2^{\kappa }+K_{1}+1}$
subsets of $[-\gamma 2^{\kappa },K_{1}]$.
So, we can choose $p_{1}''$ before \eqref{**a7}
so that we have $\# N_{A}<\epsilon _{1}q$
for all subsets $A$ of $[-\gamma 2^{\kappa },K_{1}].$

If $n\not \in N_{j}$, $n\in [0,q)$
then by \eqref{**3} and \eqref{*6x}
\begin{align}\label{*7dfw}
\#((n+\LLL_{0}(q))\cap \LLL^{j}(q))&<
K_{2}(\frac{1}{2^{\kappa }}+\rho _{1})<\\
& \nonumber
(1+\rho _{1})
(1+\rho _{1}2^{\kappa })
\gamma K_{1}(
\frac{1}{2^{\kappa }}+\rho _{1}).
\end{align}
On the other hand, by \eqref{*a10}
\begin{equation}\label{*6c}
\frac{q-\epsilon _{1}q}{K_{1}}\leq
\# \caj< \frac{q}{K_{1}}+1.
\end{equation}
For $0\leq n<q$ set $\caj(n)= \{ j\in\caj:
n\not\in N_{j}  \}$.

Consider an  $n$
such that
\begin{equation}\label{*14budd}
\# \caj(n)>q\frac{1-\epsilon _{1}}{K_{1}}
(1-\sqrt {\epsilon _{1}}).
\end{equation}
Later we show that for most $n$'s this inequality
holds.

Observe that if $x\in [-(n_{j}+K_{1})+\gamma 2^{\kappa }-1,
-n_{j}-1]\cap \LLL(q)$ then there exists $i\in \{ 0,...,\gamma 2^{\kappa }-1
  \}$ such that $x-i\in -\LLL_{0}(q)$ and
$x-i\in [-(n_{j}+K_{1}),-n_{j}-1]$, that is,
$x-i\in -\LLL_{0}^{j}(q)$ and hence $x\in \LLL^{j}(q).$
Therefore, $[-(n_{j}+K_{1})+\gamma 2^{\kappa }-1,-n_{j}-1]\cap \LLL(q)
\sse \LLL^{j}(q).$

We want to estimate
$$\#\bigg(((n+\LLL_{0}(q))\sm\LLL(q))\cap
[0,q)\bigg)
=
\#\bigg(((n+\LLL_{0}(q))\sm\LLL(q))\cap
(-q,0]\bigg).$$

By \eqref{*7dfw} for any $j\in \caj(n)$
\begin{align}\label{*15budd}
\#((n+\LLL_{0}(q))\cap& \LLL(q)\cap
[-(n_{j}+K_{1})+\gamma 2^{\kappa }-1,-n_{j}-1)<\\
\nonumber &(1+\rho _{1})(1+\rho _{1}2^{\kappa })\gamma  K_{1}
(\frac{1}{2^{\kappa }}+\rho _{1}).
\end{align}

Put $$T_{n}= \{ t\in \Z:
t\in (-q,0]\sm \bigcup_{j\in \caj(n)}
[-(n_{j}+K_{1})+\gamma 2^{\kappa }-1,-n_{j}-1]  \}.$$
It is clear that
\begin{align}\label{*16budd}
&
(n+\LLL_{0}(q))\cap \LLL(q)
\cap (-q,0]\sse\\ \nonumber
&T_{n}\cup
\bigcup_{j\in\caj(n)}
(n+\LLL_{0}(q))\cap \LLL(q)\cap
[-(n_{j}+K_{1})+\gamma 2^{\kappa }-1,-n_{j}-1]
.
\end{align}

We need to estimate $\# T_{n}$.
Since the intervals $[-(n_{j}+K_{1}),-n_{j}-1]$
are disjoint and, with the possible exception of
the one with the largest index, are subsets of
$(-q,0]$ we have by using \eqref{**3} and
\eqref{*14budd}
$$\# T_{n}< q-\# \caj(n) (K_{1}-\gamma 2^{\kappa })<
q-\# \caj(n) (1-\epsilon _{1})K_{1}<
$$
$$q-q\frac{(1-\epsilon _{1})^2}{K_{1}}(1-\sqrt{\epsilon _{1}})K_{1}=
q(1-(1-\epsilon _{1})^2(1-\sqrt {\epsilon _{1}}))<q
c_{1}(\epsilon _{1}),$$
where $c_{1}(\epsilon _{1})\to 0$ as $\epsilon _{1}\to 0.$
Now we use this, \eqref{*15budd} and \eqref{*16budd}
to estimate
$$\# \bigg (
(n+\LLL_{0}(q))\cap \LLL(q)
\cap (-q,0]\bigg)<$$ $$
\#\caj(n)\cdot (1+\rho _{1})(1+\rho _{1}2^{\kappa })
\gamma K_{1}(\frac{1}{2^{\kappa }}+\rho _{1})+q c_{1}(\epsilon _{1})<
$$
(using \eqref{*6c})
$$(\frac{q}{K_{1}}+1)(1+\rho _{1})
(1+\rho _{1}2^{\kappa })\gamma K_{1}(\frac{1}{2^{\kappa }}+\rho _{1})
+q c_{1}(\epsilon _{1})=$$ $$
(1+\frac{K_{1}}q)(1+\rho _{1})(1+\rho _{1}2^{\kappa })^{2}
\gamma q \frac{1}{2^{\kappa }}+q c_{1}(\epsilon _{1})<$$
(using \eqref{*5budd} and $q>p_{1}'>K_{1}/\epsilon _{1}$)
$$(1+\epsilon _{1})(1+\rho _{1})(1+\rho _{1}2^{\kappa })^{2}
\gamma q \frac{1}{2^{\kappa }}+q c_{1}(\epsilon _{1})=$$
$$ \bigg((1+\epsilon _{1})(1+\rho _{1})(1+\rho _{1}2^{\kappa })^{2}
+\frac{2^{\kappa }}{\gamma } c_{1}(\epsilon _{1})\bigg)
\gamma q \frac{1}{2^{\kappa }}=
c_{2}(\epsilon _{1},\rho _{1})\gamma  q \frac{1}{2^{\kappa }},$$
where $c_{2}(\epsilon _{1},\rho _{1})\to 1$ as $\epsilon _{1},\rho _{1}
\to 0.$

We can choose $\epsilon _{1},
\rho _{1}>0$ so that
$$c_{2}(\epsilon _{1},\rho _{1})
<\frac1{(1-\rho _{0})}.
$$
This by \eqref{*5dfw} implies
$$c_{2}(\epsilon _{1},\rho _{1})\gamma q \frac{1}{2^{\kappa }}
<\frac1{(1-\rho _{0})^2}\gamma \#
(\LLL_{0}(q)\cap [0,q)).$$

By \eqref{*6budd} we obtain
\begin{align}\label{*17budd}
&\#\bigg(
((n+\LLL_{0}(q))\sm \LLL(q))\cap [0,q)
\bigg)>\\
&\#(\LLL_{0}(q)\cap [0,q))(1-\frac1{(1-\rho _{0})^{2}}\gamma )>
\nonumber
\#(\LLL_{0}(q)\cap [0,q))(1-\rho )(1-\gamma ).
\end{align}
To prove \eqref{*a9} we need to show that there are sufficiently
many $n$'s which satisfy the above inequality.

 Let $b$ be the number of $n$'s for which
\begin{equation}\label{*5}\# \caj(n)\leq
q\frac{1-\epsilon_{1}}
{K_{1}}(1-\sqrt{\epsilon_{1}}).
\end{equation}
If we can show that $b<\epsilon q$, then we have finished the
proof of Lemma \ref{probabilitylemma} since if \eqref{*14budd}
holds for an $n$ then we have \eqref{*17budd}.
If $n$ satisfies \eqref{*5} then by the definition
of $\caj(n)$
from
\eqref{*6c} we infer that
 $n\in N_j$ for at least $\sqrt{\epsilon _{1}}
\frac{1-\epsilon _{1}}
{K_{1}}q$ many $j$'s. Hence, using \eqref{**5} and \eqref{*6c}
$$b\sqrt{\epsilon _{1}}
\frac{1-\epsilon _{1}}{K_{1}}q\leq
\sum_{j\in \caj}\# N_{j}\leq
\epsilon _{1}q(\frac{q}{K_{1}}+1)< \frac{2\epsilon _{1}q^{2}}
{K_{1}}
$$
which implies
$$b\leq
\frac{2\sqrt{\epsilon _{1}}}{1-\epsilon _{1}}q<\epsilon q.$$
\end{proof}

\section{Periodic rearrangements}\label{perr}

We need a lemma concerning the fact  that  one can make a periodic
perturbation of certain given periodic sets so that averages of
the characteristic function of the perturbed set taken along squares
is close  to the average
measure of the original set.

Assume $F\sse \R$ is periodic by $\tau '\in \N$
and if $x\in F$ then $[\lf x\rf,\lf x\rf +1)\sse F.$
Given a
natural
number $\tau>\tau '$, the
{\it $\tau$-periodic rearrangement}
of $F$ is denoted by $F^{\tau }$ and it is periodic by $\tau ,$
and $F^{\tau }\cap [0,\tau )=[0,\lf \tau /\tau '\rf \cdot \tau ')
\cap F.$

For the proof of the next lemma we recall that by the P\'olya-Vinogradov Theorem
(see for example p. 324 of \cite{[IK]}) for any $n\in \Z$, $l\in\N$ 
and odd prime $\tttt$, 
we have
$\ds \sum_{j=n}^{n+l-1}\bigg(\frac{j}{\tttt}\bigg)_{\call}\leq 
6\sqrt \tttt \log\tttt .$

\begin{lemma}\label{taurearrangement}
Suppose  $\tau \in \N$, $F \subset \R$ periodic by $\tau '$ and
$\rho >0$ are given. There exists
$M_{\rho }$ such that if $\tau >M_{\rho }$
is a prime number then for any
$n\in \N$ and
$x\in \R$,
if $\tau| m$ then
\begin{equation}\label{*d1}
\frac{1}{m}\sum_{k=n}^{n+m-1}
\chi_{F^{\tau }}(x+k^{2})\geq (1-\rho )\olll(F).
\end{equation}
\end{lemma}

\begin{proof} For each $a \in \{1,...,\tau'\}$, let $\ds F_a =
  \cup_{n\equiv a \ 
    (mod\  \tau')} [n, n+1)$. Since $F$ is the disjoint union of
    finitely many of the sets
  $F_a$, it suffices to prove the lemma for some fixed $F_a$.
Also, note that it is enough to show that there is some $M_{\rho }$ such
  that if $\tau >M_{\rho }$
is a prime number, then for any
$n\in \N$ and
$x\in \R$,

\begin{equation}
\frac{1}{\tau}\sum_{k=n}^{n+\tau-1}
\chi_{F_a^{\tau }}(x+k^{2}) \geq (1-\rho )\olll(F_a) = \frac{1-\rho}{\tau'}.
\end{equation}

Put $l=\lf\tttt/\tttt' \rf$.

Note $\chi_{F_a^\tau}(x+k^2) = 1$ for an integer $x$
 if and only if there is some
$j\in\{0,...,l-1\}$ such that $x+k^2
\equiv a+j\tau' \ (mod\ \tau).$

Thus,  using the Legendre symbol $\ds \bigg (\frac{k}{n}\bigg)_\call$, we have

$$
\frac{1}{\tau} \#\{k\in \{0,...,\tau -1\}:k^2 \equiv a+j\tau'-x, \ (mod \  \tau)
\text{ for a }j\in \{ 0,...,l-1 \}\}
=
$$
$$
\frac{1}{\tau}\sum_{j=0}^{l-1}\#\{k\in \{0,...,\tau -1\}:k^2 \equiv
a+j\tau'-x, \  (mod \ 
\tau)\}=
$$
$$
\frac{1}{\tau}\sum_{j=0}^{l-1}\bigg(1+\bigg(\frac{a+j\tau'-x}{\tau}\bigg)
_\call\bigg) =
\frac{l}{\tau} + \frac{1}{\tau}\sum_{j=0}^{l-1}\bigg(
\frac{a+j\tau'-x}{\tau}\bigg)_\call=
$$
$$
 \frac{\lf{\tau}/{\tau'}\rf}{\tau} + \frac{1}{\tau}\sum_{j=0}^{l-1}
 \bigg(\frac{a+j\tau'-x}{\tau}\bigg)_\call
$$
Now as $\tau \rightarrow \infty$, $
\frac{\lf{\tau}/{\tau'}\rf}{\tau} \rightarrow \frac{1}{\tau'}.$ Setting $\ds S
=\sum_{j=0}^{l-1}\bigg (\frac{a+j\tau'-x}{\tau}\bigg)_{\call}$ we only need to show
$S/\tau \rightarrow 0$, as $\tau \rightarrow \infty$.

We argue this as follows. Since $\tau'$ is a prime with $\tau' < \tau$,
choose $b$ such that $b\tau' \equiv a-x,\  (mod \ \tau)$ and set $\tau^* =
(\tau')^{-1},\   (mod \ \tau)$ so that $b \equiv (a-x)\cdot\tau^*,\   (mod \
\tau)$. Since the Legendre symbol is a character, 
$$\bigg(\frac{(a+j\tau'-x)\cdot\tau^*}{\tau}\bigg)_\call =  
\bigg(\frac{a+j\tau'-x}{\tau}\bigg)_\call
\cdot\bigg(\frac{\tau^*}{\tau}\bigg)_\call.$$ Also, we have
$|\big(\frac{\tau^*}{\tau}\big)_\call | = 1.$ Since 
$$\bigg(\frac{(a+j\tau'-x)\cdot \tau^*}{\tau}\bigg)_\call = 
\bigg(\frac{(a-x)\tau^*+j\tau'\tau^*}{\tau}\bigg)_\call 
=\bigg(\frac{b+j}{\tau}\bigg)_\call,$$ we have
by the P\'olya-Vinogradov inequality 
$$ 
|S| = \bigg|\sum_{j=0}^{l-1}\bigg(\frac{b+j}{\tau}\bigg)_\call \bigg| 
\leq 6\sqrt\tau\log\tau \text{ and }\frac{|S|}{\tttt}\to 0
\text{ as }\tttt\to\oo.
$$
This completes the proof of Lemma 3.

\end{proof}

\section{$K-M$  families}\label{secksystem}

\begin{definition}\label{conmkilenc}
For a positive integer $M$
we say that a periodic function
or a ``random variable",
$X:\R\to \R$ is {\it conditionally $M$$-$$0.99$ distributed}
 on the
set $\LLL$, which is
 periodic by the same period, if
$X(x)\in  \{ 0,0.99,0.99\cdot \frac{1}2,...,
0.99\cdot 2^{-M+1}  \}$,
 and
$\olll( \{ x\in\LLL:X(x)=0.99\cdot 2^{-l}  \})
=0.99 \cdot 2^{-M+l-1}\olll(\LLL)$ for $l=0,...,M-1.$
(We regard $\R$
as being
 periodic by $1$ with $\olll(\R)=1$
and if $\LLL=\R$ then we just simply say
that $X$ is $M$$-$$0.99$-distributed.)
By an obvious adjustment this definition will also be used
for random variables $\oXX$ defined on $[0,1)$ equipped
with the Lebesgue measure $\lambda .$
If we have two ``random variables" $X_{1}$ and $X_{2}$
both conditionally $M$$-$$0.99$ distributed on $\LLL$
then they are called pairwise independent (on $\LLL$) if
for any $y_{1},y_{2}\in \R$
\begin{align}\label{indep}
\olll \{ x\in\LLL &:X_{1}(x)=y_{1}\text{ and } X_{2}(x)=y_{2} \}
\olll(\LLL)=\\&
\nonumber \olll( \{ x\in\LLL:X_{1}(x)=y_{1}  \})\olll( \{ x
\in\LLL:
X_{2}(x)=y_{2}  \})\\
\intertext{ or, equivalently,}
\nonumber
\olll \{ x\in\LLL&:X_{1}(x)=y_{1}\text{ and } X_{2}(x)=y_{2} \}
/\olll(\LLL)=\\ \nonumber
& (\olll( \{ x\in\LLL:X_{1}(x)=y_{1}  \})/\olll(\LLL))
(\olll( \{ x\in\LLL:
X_{2}(x)=y_{2}  \})/\olll(\LLL)).
\end{align}
\end{definition}

If we say that $X_{1}$ and $X_{2}$ are pairwise independent,
without specifying $\LLL$ then we mean $\LLL=\R$.

We will use the following simple properties.
Assume $\LLL_{1}$ and $\LLL_{2}$ are two disjoint sets
with a common period.
If $X_{1}$ and $X_{2}$
are conditionally $M$$-$$0.99$ distributed on $\LLL_{1}$ and
on $\LLL_{2}$, then $X_{1}$ (and similarly $X_{2}$) is
conditionally $M$$-$$0.99$ distributed on $\LLL_{1}\cup
\LLL_{2}$.
If, in addition $X_{1}$ and $X_{2}$ are pairwise independent
on each $\LLL_{1}$ and $\LLL_{2}$, then $X_{1}$ and $X_{2}$ are
pairwise independent on $\LLL_{1}\cup
\LLL_{2}$.
We note the
last property depends on $X_{1}$ and $X_{2}$
having the {\it same distribution} on $\LLL_{1}$ and $\LLL_{2}$.
Similar properties hold if we have finitely many functions
$X_{1},...,X_{K}$ with the same conditional distribution.

For our argument a wide range of independent, identically distributed
uniformly bounded ``random variables" with expectations
bounded from below by constant times $M2^{-M}$ could be used. 
However, as the remark above shows we need identically distributed
random variables and 
out of  the many possible
choices we picked the $M-0.99$ distributed ones. For a motivation
for this choice see \cite{BM}.

We say that $X:\R\to \R$ is $M-0.99$ super distributed if
\begin{equation}\label{*supdist1}
X(x)\in  \{ 0,0.99,0.99\cdot 2^{-1},...,
0.99\cdot 2^{-M+1}  \},
\end{equation}
 and
$$\olll( \{ x\in\R:X(x)=0.99\cdot 2^{-l}  \})
\geq 0.99 \cdot 2^{-M+l-1} \text{ for }l=0,...,M-1.$$
We need the following lemma:
\begin{lemma}\label{*super}
Suppose $\tttt\in \N$, $X_{1},...,X_{K}:\R\to [0,\oo)$
are $M-0.99$ distributed, $\tttt$ periodic and $X_{K+1}'$
is $M-0.99$ super distributed $\tttt$ periodic and $X_{K+1}'$
is pairwise independent from $X_{h}$ for all $h=1,...,k.$
Then we can choose $0\leq X_{K+1}\leq X_{K+1}'$ such that
$X_{K+1}$ is $M-0.99$ distributed and pairwise 
independent from $X_{h}$ for all $h=1,...,k.$
\end{lemma}
\begin{proof}
Set $$\TTT_{K+1,l}'\defeq \{ x\in \R:X_{K+1}'(x)=0.99\cdot 2^{-l} \}.$$
Then $\olll(\TTT_{K+1,l}')\geq 0.99\cdot 2^{-M+l-1}$ and
$$1\geq c_{l}\defeq \frac{0.99 \cdot 2^{-M+l-1}}{\olll(\TTT_{K+1,l}')}.$$
We also set $\cas_{K}=\{ 0,0.99,...,0.99\cdot 2^{-M+1} \}^{K}$
and for $(y_{1},...,y_{K})\in \cas_{K}$ set
$$\TTT_{K+1,l}'(y_{1},...,y_{K})\defeq \{ x\in \TTT_{K+1,l}': X_{h}(x)=y_{h},
\ h=1,...,K \}.$$
Clearly, $\TTT_{K+1,l}'(y_{1},...,y_{K})$ is $\tttt$ periodic and for
$(y_{1},...,y_{K})\not =(y_{1}',...,y_{K}')\in \cas_{K}$
the sets $\TTT_{K+1,l}'(y_{1},...,y_{K})$ and $\TTT_{K+1,l}'(y_{1}',...,y_{K}')$
are disjoint. For all $(y_{1},...,y_{K})\in \cas_{K}$ choose a Borel measurable
$\TTT_{K+1,l}(y_{1},...,y_{K})\sse \TTT_{K+1,l}'(y_{1},...,y_{K})$
such that $$\olll(\TTT_{K+1,l}(y_{1},...,y_{K}))=
c_{l}\olll (\TTT_{K+1,l}'(y_{1},...,y_{K}))$$
and $\TTT_{K+1,l}(y_{1},...,y_{K})$ is periodic by $\tttt$.

For $x\in \TTT_{K+1,l}(y_{1},...,y_{K})$ set $X_{K+1}(x)=X_{K+1}'(x)=
0.99\cdot 2^{-l}$ and for $x\in \TTT_{K+1,l}'(y_{1},...,y_{K})\sm
\TTT_{K+1,l}(y_{1},...,y_{K})$ set $X_{K+1}(x)=0.$
Do this for all $(y_{1},...,y_{K})\in \cas_{K}$ and for all $l=0,...,M-1.$
Finally, for those $x$'s for which $X_{K+1}'(x)=0$ set $X_{K+1}(x)=0$.
Then $0\leq X_{K+1}\leq X_{K+1}'.$ 

Suppose $l\in \{ 0,...,M-1 \}$ is fixed. Then
$$\olll(\{ x: X_{K+1}(x)=0.99\cdot 2^{-l} \})=
\sum_{(y_{1},...,y_{K})\in \cas_{K}}\olll(\TTT_{K+1,l}(y_{1},...,y_{K}))=
$$
$$\sum_{(y_{1},...,y_{K})\in \cas_{K}}c_{l}\olll(\TTT_{K+1,l}'(y_{1},...,y_{K}))=
c_{l}\olll(\TTT_{K+1,l}')=0.99\cdot 2^{-M+l-1}.$$
This and \eqref{*supdist1} also implies that 
\begin{equation}\label{*komplest}
\olll(\{ x : X_{K+1}(x)=0 \})=1-\sum_{l=0}^{M-1}
0.99\cdot 2^{-M+l-1}.
\end{equation}

Suppose $\oyy_{h}\in \{ 0,0.99,0.99\cdot 2^{-1},...,
0.99\cdot 2^{-M+1}  \}$ is fixed and denote by $\cas_{K,\oyy_{h}}$
the set of those $(y_{1},...,y_{K})\in \cas_{K}$ for which $y_{h}=\oyy_{h}.$
Then by the pairwise independence of $X_{h}$ and $X_{K+1}'$ on $\R$
we have
\begin{equation}\label{*indep171}
\olll(\{ x: X_{h}(x)=\oyy_{h} \text{ and }
X_{K+1}'(x)=0.99\cdot 2^{-l} \})=
\end{equation}
$$\olll(\{ x:X_{h}(x)=\oyy_{h} \})\olll(\{ x : X_{K+1}'(x)=0.99\cdot 2^{-l}\})=
$$ $$
\olll(\{ x:X_{h}(x)=\oyy_{h} \})\olll(\TTT_{K+1,l}').$$
On the other hand,
\begin{equation}\label{*indep172}
\sum_{(y_{1},...,y_{K})\in\cas_{K,\oyy_{h}}}\olll(\TTT_{K+1,l}'(y_{1},...,y_{K}))=
\end{equation}
$$\olll(\{ x:X_{h}(x)=\oyy_{h} \text{ and } X_{K+1}'(x)=0.99\cdot 2^{-l}\}).$$
Since
$$
c_{l}\sum_{(y_{1},...,y_{K})\in\cas_{K,\oyy_{h}}}\olll(\TTT_{K+1,l}'(y_{1},...,y_{K}))=
\sum_{(y_{1},...,y_{K})\in\cas_{K,\oyy_{h}}}\olll(\TTT_{K+1,l}(y_{1},...,y_{K}))=$$
$$
\olll(\{ x:X_{h}(x)=\oyy_{h} \text{ and } X_{K+1}(x)=0.99\cdot 2^{-l}\}),$$
if we multiply \eqref{*indep171} and \eqref{*indep172} by $c_{l}$
we obtain
$$\olll(\{ x:X_{h}(x)=\oyy_{h} \text{ and } X_{K+1}(x)=0.99\cdot 2^{-l}\})=$$
$$\olll(\{ x:X_{h}(x)=\oyy_{h} \})\cdot c_{l}\cdot \olll(\TTT_{K+1,l}')=
$$
$$\olll(\{ x:X_{h}(x)=\oyy_{h} \})\cdot 0.99\cdot 2^{-M+l-1}=$$
$$\olll(\{ x:X_{h}(x)=\oyy_{h} \})\olll(\{ x : X_{K+1}(x)=0.99\cdot 2^{-l}\}).$$
By \eqref{*supdist1} and \eqref{*komplest} we also have
$$\olll(\{ x:X_{h}(x)=\oyy_{h} \text{ and } X_{K+1}(x)=0\})=$$
$$\olll(\{ x:X_{h}(x)=\oyy_{h} \})-\sum_{l=0}^{M-1}
\olll(\{ x:X_{h}(x)=\oyy_{h} \text{ and } X_{K+1}(x)=0.99\cdot 2^{-l}\})=$$
$$\olll(\{ x:X_{h}(x)=\oyy_{h} \})(1-\sum_{l=0}^{M-1}0.99\cdot 2^{-M+l-1})=$$
$$\olll(\{ x:X_{h}(x)=\oyy_{h} \})\olll(\{ x : X_{K+1}(x)=0\}).$$
This completes the proof of the fact that $X_{K+1}$
is pairwise independent from $X_{h}$ for all $h=1,...,K.$
\end{proof}

\begin{definition}\label{nagykompl}
We say that a set $\capp\sse \N$ has {\it sufficiently large
complement} if
there are infinitely many primes
relatively prime to any number in $\capp.$
\end{definition}

Sometimes we need to work with the ``real'' squares modulo $q$:

\begin{definition}\label{oLLLdfw}
Assume $q=p_{1}\cdot \cdot \cdot p_{\kappa }$, where
$p_{1}<...<p_{\kappa }$ are odd primes.
Set $$\LLL_{0}'(q)= \{ n\in \LLL_{0}(q): p_{j}\not|n, \text{ for all
}j=1,...,\kappa  \}.$$
If $n\in \LLL_{0}'(q)$, then there are $2^{\kappa }$ many solutions
of $x^{2}\equiv n$ $mod$ $q$, also observe that for fixed $\kappa $
\begin{equation}\label{*9b}
\text{ if }
p_{1}\to \oo, \text{ then }\frac{\#(\LLL_{0}(q)\cap [0,q))}
{\#(\LLL_{0}'(q)\cap [0,q))}\to 1.
\end{equation}
Given $\gamma \in (0,1)$ we also put
$$\LLL'_{\gamma  }(q)
=-\LLL'_{0}(q)+ \{ j\in\Z:0\leq j <\gamma 2^{\kappa }  \},$$
$$\oLLL_{\gamma  }'(q)=\LLL_{\gamma  }'(q)+[0,1)=-\LLL_{0}'(q)+
 \{ x:0\leq x <\gamma 2^{\kappa }  \}.$$
In the sequel often if $\gamma$ is fixed
we will suppress the dependence on $\gamma $
by writing $\LLL'(q)$ and $\oLLL'(q)$ instead of
$\LLL'_{\gamma  }(q)$ and $\oLLL'_{\gamma  }(q)$,
respectively.
 To help to memorize our notation of these sets, ``$\LLL'$" 
means that the set ``$\LLL$" is built as ``$\LLL$" but instead
of $\LLL_0$ we use $\LLL_{0}'$ in our construction. We keep our earlier
convention as well and hence ``$\oLLL'$" is the set obtained 
from ``$\LLL'$" by adding $[0,1)$.
\end{definition}

\begin{definition}\label{kmfamily}
Suppose $K,M\in \N$, $\LLL\sse \R$ is periodic by $\tqq.$
There is a parameter $\gamma' $ associated to $\LLL$.
(If $\LLL=\R$ then $\tqq=\gamma' =1$. Otherwise
one should think of $\Lambda=\oLLL'_{\gamma '}
(\tqq)$ and $\gamma '$ is the parameter
used
in the definition
of $\Lambda$.)
In the sequel we assume that $\capp\sse \N$ has
sufficiently large complement.
A $K-M$ {\it family living on}
$\LLL$ with input parameters $\delta >0,$ $\OOO$, $\GGG>1$,
$A\in \N,$  $\capp$  with output objects
$\tau $, $f_{h}$,  $X_{h}$ ($h=1,...,K$); $E_{\delta }$,
$\omega (x)$, $\alpha (x)$ and $\tau (x)$ is a system
satisfying:

\begin{enumerate}
\item There exist a period $\tau $,
functions
 $f_{h}:
\R\to [0,\oo)$,
pairwise independent, conditionally
$M$$-$$0.99$-distributed on $\LLL$ ``random" variables
$X_{h}:\R\to\R$, for $h=1,...,K,$
and a set $E_{\delta }$
such that
 all these objects are periodic by $\tau$ where
$\tau$ is an integer multiple of $\tqq$.
\item
We have $\olll(E_{\delta })<\delta$.
For all $x\not\in E_{\delta }$,
there exist
 $\omega (x)>\alpha (x)>A,$ $\tau (x)< \tau $
such that
$\omega ^{2}(x)<\tau ,$
$\frac{\omega (x)}{\alpha (x)}>\OOO \cdot  \tau (x)$;
moreover
if $\alpha (x)\leq n<n+m\leq \omega (x)$
and
$\tau (x)|m$, then for all $h=1,...,K,$
\begin{equation}\label{*12dfw}
\frac{1}m\sum_{k=n}^{n+m-1}f_{h}(x+k^{2})>X_{h}(x).
\end{equation}
\item
For all $p\in\capp$, $(\tau (x),p)=1,$ $(\tau ,p)=1.$
\item
For all $x\in
\Lambda\sm E_{\delta }$, for all $h\in  \{ 1,...,K  \}$
\begin{equation}\label{**9b}
f_{h}(x+j+\tau (x))=f_{h}(x+j)
\end{equation}
 whenever
$\alpha ^{2}(x)\leq j<j+\tau (x)\leq \omega ^{2}(x).$
\item Finally, for $h=1,...,K$
\begin{equation}\label{*17iest}
\frac{1}\tau \int_{0}^{\tau }f_{h}=\ointt f_{h}<\GGG\cdot
\gamma' \cdot
2^{-M+1}.
\end{equation}
\end{enumerate}
\end{definition}

\begin{remark}
The input parameters in the above definition should
be regarded as something given in advance while the
output objects are defined and constructed later.
The most important property is \eqref{*12dfw},
while the numerous other technical properties are
needed in order to verify by mathematical induction
the existence of $K-M$ families.

If $x$ is not in the exceptional set $E_{\delta }$,
then
\eqref{*12dfw} says that
the average of $f_{h}$ 
taken along the squares of a run
of integers staying in the
window $[\alpha (x),\omega (x)]$ dominates $X_{h}(x)$,
provided that the length of the run is a multiple
of $\tttt(x)$.
 In \eqref{**9b} we claim that
these functions appear to be periodic in the window
$[\alpha ^{2}(x),\beta ^{2}(x)],$
that is, when squares stay in the window
$[\alpha(x),\beta (x)].$
\end{remark}

\begin{lemma}\label{lksystem}
Let $M\in \N$, $\LLL=\R$ (this implies
$\tqq=\gamma' =1$). Then for each positive integer
$K$ and parameters
 $\delta >0,$ $\OOO,\GGG>1$, $A\in \N,$
and $\capp \sse \N$
such that
$\capp$ has sufficiently large complement
there exist a $K-M$ family living on $\R$
with these parameters.
\end{lemma}

\subsection{Putting $K-M$ families on quadratic residue classes}\label{41}

The  proof of Lemma \ref{lksystem} is quite involved.
It will be done by induction on $K$.
We will build a $(K+1)-M$ family for a given set of
input parameters, provided we know
the existence of $K-M$ families for all possible input parameters.
To carry out this  induction step we need to verify that
a generalized version of Lemma \ref{lksystem} holds.

So
 we assume that $K-M$ families
on $\R$
exist for a fixed
$K\in \N$ for all possible parameter values.

We will use the following lemma about {\it ``putting a
$K-M$-family on a residue class".}
We assume that $\ocap$ is a set of natural numbers
with sufficiently large complement
and we have a number $\tqq$ such that
\begin{equation}\label{50315a}
(\tqq  ,p)=1 \text{ for all } p\in\ocap,\text{  }
\tqq  =p_{0,1}\cdot \cdot \cdot p_{0,\kappa },
\end{equation}
and
$p_{0,1}<...<p_{0,\kappa }$
are odd primes.

We also assume that a constant $\gamma =2^{-c_{\gamma }}$,
the
so called ``leakage constant" is given with $c_{\gamma }\in \N$
and $\kappa>c_{\gamma}$.
This $\gamma $
is
 used in the definition of $\oLLL'(\tqq  )=\oLLL'_{\gamma  }(\tqq)$.

\begin{lemma}\label{lputres}
Let $M\in \N$ be given and suppose for some $K\in \N$
that
$K-M$ families exist
on $\R$
for all possible parameter
values.
Suppose that
$\ocap$,
$\tqq$ and the parameter $\gamma$
associated to $\oLLL'(\tqq)$ satisfy the above assumptions.
In addition, let
$\delta >0,$ $\OOO,\GGG>1$, and $A\in\N$
be given.
Then
for the above input parameters
there exists a $K-M$ family living on
$\oLLL'(\tqq  )$
with output objects $\ottt$, $\off_{h}$, 
$\oXX_{h}$ ($h=1,...,K$); $\oEE_{\delta }$, $\oooo(x)$,
$\oaaa(x)$ and $\ottt(x)$. Moreover, 
$\ottt=\tttt \tqq$ with a suitable $\tttt\in\N$ and
if $\oXXX(\tqq)=
\cup_{j\in \Z} [j\tqq ,j\tqq +\gamma 2^{\kappa })$ then
$\off_{h}(x)=0$ for $x\not\in \oXXX(\tqq)$ and $h=1,...,K.$
\end{lemma}

\begin{proof}
Using $\capp=\ocap\cup \{ \tqq   \}$,
choose a $K-M$ family living on $\R$
with input parameters
$\delta ,$
$\OOO'=\OOO \tqq ,$
$\GGG,$ $A.$
This $K-M$ family provides us
with
$\tau $, $f_{h}$,
 $X_{h}$,
$E_{\delta },$ $\omega (x),$ $\alpha (x)$, and $\tau (x)$,
 satisfying (i)-(v) of Definition \ref{kmfamily}.
Especially,
\begin{equation}\label{*14dfww}
\text{for all } p\in \ocap\cup  \{ \tqq   \}
\text{ we have }(\tau (x),p)=1\text{ and }(\tau ,p)=1.
\end{equation}
We construct a new $K$-system, marked by overlines,
which lives on $\oLLL'(\tqq)$ and is periodic by
$\ottt=\tau \tqq$.

Set $\off_{h}(x)=f_{h}(x)\tqq /2^{\kappa }$ if $x\in
\oXXX(\tqq)=\cup_{j\in \Z} [j\tqq ,j\tqq +\gamma 2^{\kappa })$,
otherwise put $\off_{h}(x)=0.$
Then $\off_{h}$ is periodic by $\tau \tqq $.

Next we define $\oXX_{h}(x)$ so that
$\oXX_{h}(x)=X_{h}(x)$
for $x\in\oLLL'(\tqq )$, and
otherwise let $\oXX_{h}(x)=0.$
Clearly, $\oXX_{h}$ is periodic by $\tau \tqq$ and is
supported on $\oLLL'(\tqq)$.

Next we check the
distribution of $\oXX_{h}(x)|_{\oLLL'(\tqq )}.$
We know that $X_{h}(x+\tau )=X_{h}(x)$.
From \eqref{*14dfww},  $(\tau , \tqq)=1$
and thus 
the numbers
$j\tau $, $j=0,...,\tqq -1$ cover all residue classes modulo
$\tqq .$
Now we can compute
$$\olll ( \{ x\in \oLLL'(\tqq):
\oXX_{h}(x)=0.99\cdot 2^{-l}  \})=$$
$$\olll ( \{ x\in \R:
\oXX_{h}(x)=0.99\cdot 2^{-l}  \})=$$
$$\frac1{\tau \tqq }\lambda ( \{ x\in [0,\tau \tqq ):
\oXX_{h}(x)=0.99\cdot 2^{-l}  \})=$$
$$\frac1{\tau \tqq }
\sum_{j=0}^{\tqq -1}
\lambda ( \{ x\in [j\tau ,(j+1)\tau ):
\oXX_{h}(x)=0.99\cdot 2^{-l}  \})=$$
$$
\frac1{\tau \tqq }
\sum_{j=0}^{\tqq -1}
\lambda ( \{ x\in [0,\tau ):
\oXX_{h}(x+j\tau )=0.99\cdot 2^{-l}  \})=
$$
$$
\frac1{\tau \tqq }
\# (\LLL'(\tqq )\cap [0,\tqq))\cdot
\lambda ( \{ x\in [0,\tau ):
X_{h}(x)=0.99\cdot 2^{-l}  \})=$$
(using that $\#(\LLL'(\tqq) \cap [0, \tqq))/\tqq =\olll (\oLLL'(\tqq ))$)
$$
\olll(\oLLL'(\tqq))
\cdot \frac1{\tau }\lambda( \{ x\in[0,\tau ):X_{h}(x)=0.99\cdot 2^{-l}  \})=
$$
$$\olll(\oLLL'(\tqq))
\cdot \olll( \{ x\in\R:X_{h}(x)=0.99\cdot 2^{-l}  \})=
\olll(\oLLL'(\tqq ))\cdot 0.99\cdot 2^{-M+l-1}.$$
Thus the ``conditional distribution"
of $\oXX_{h}$ on $\oLLL'(\tqq )$
is $M$$-$$0.99.$

Set $\Upsilon = \{ 0, 0.99, 0.99\cdot 2^{-1},...,0.99\cdot 2^{-M+1}  \}$
and $\Upsilon _{+}=\Upsilon \sm  \{ 0  \}.$
Next we show that the functions $\oXX_{h}$ are pairwise independent
on $\oLLL'(\tqq)$.
Suppose $h_{1}\not=h_{2}$.
First assume $y_{1},y_{2}\in \Upsilon _{+}$.
Then the above argument shows
\begin{equation}\label{*50314a}
\olll( \{ x\in \oLLL'(\tqq):
\oXX_{h_{j}}(x)=y_{j}  \})=
\olll( \{ x\in \R:
X_{h_{j}}(x)=y_{j}  \})
\olll(\oLLL'(\tqq))
\end{equation}
for $j=1,2.$
A similar
 argument
 shows
\begin{align}\label{*50314b}
&\olll( \{ x\in \oLLL'(\tqq):
\oXX_{h_{1}}(x)=y_{1}
\text{ and }
\oXX_{h_{2}}(x)=y_{2} \})=\\
& \nonumber
\olll( \{ x\in \R:
X_{h_{1}}(x)=y_{1}
\text{ and }
X_{h_{2}}(x)=y_{2}\})
\olll(\oLLL'(\tqq)).
\end{align}

The range of $\oXX_{h_{j}}$ and $X_{h_{j}}$ equals
$\Upsilon =\Upsilon _{+}\cup  \{ 0  \}$ and
\eqref{*50314a} holds for all $y_{j}\in \Upsilon _{+}$.
Therefore,
\begin{equation}\label{*50314c}
\olll( \{ x\in \oLLL'(\tqq):
\oXX_{h_{j}}(x)=0  \})=
\olll( \{ x\in \R:
X_{h_{j}}(x)=0  \})
\olll(\oLLL'(\tqq))
\end{equation}
should also hold for $j=1,2.$

Recalling that $X_{h_{j}}$ are $M-0.99$ distributed
and pairwise independent
on $\R$,
using for a fixed $y_{2}\in \Upsilon _{+}$, \eqref{*50314b}
for all $y_{1}\in \Upsilon _{+}$, and using \eqref{*50314a}
with $j=2$ one can deduce
\begin{align}\label{*50314d}
&\olll( \{ x\in \oLLL'(\tqq):
\oXX_{h_{1}}(x)=0
\text{ and }
\oXX_{h_{2}}(x)=y_{2} \})=\\
& \nonumber
\olll(\{x\in\oLLL'(\tqq): \oXX_{h_{2}}(x)=y_2\})-
\\
& \nonumber
\qquad
\olll( \{ x\in \oLLL'(\tqq):
\oXX_{h_{1}}(x)\in \Upsilon _{+}
\text{ and }
\oXX_{h_{2}}(x)=y_{2} \})=
\\
& \nonumber
\olll( \{ x\in \R:
X_{h_{1}}(x)=0
\text{ and }
X_{h_{2}}(x)=y_{2}\})
\olll(\oLLL'(\tqq)).
\end{align}
Similarly, one can see
that for any $y_{1}\in \Upsilon _{+}$
\begin{align}\label{*50314e}
&\olll( \{ x\in \oLLL'(\tqq):
\oXX_{h_{1}}(x)=y_{1}
\text{ and }
\oXX_{h_{1}}(x)=0 \})=\\
& \nonumber
\olll( \{ x\in \R:
X_{h_{1}}(x)=y_{1}
\text{ and }
X_{h_{2}}(x)=0\})
\olll(\oLLL'(\tqq)).
\end{align}
Recalling that $X_{h_{2}}$ is $M-0.99$-distributed
on $\R$ and
using \eqref{*50314d} for all $y_{2}\in \Upsilon _{+}$
and \eqref{*50314c} with $j=1$ one can see that
\begin{align}\label{*50314f}
&\olll( \{ x\in \oLLL'(\tqq):
\oXX_{h_{1}}(x)=0
\text{ and }
\oXX_{h_{2}}(x)=0 \})=\\
& \nonumber
\olll(\{x\in\oLLL'(\tqq): \oXX_{h_{1}}(x)=0\})-
\\
& \nonumber
\qquad
\olll( \{ x\in \oLLL'(\tqq):
\oXX_{h_{1}}(x)=0
\text{ and }
\oXX_{h_{2}}(x)\in \Upsilon _{+} \})=
\\
& \nonumber
\olll( \{ x\in \R:
X_{h_{1}}(x)=0
\text{ and }
X_{h_{2}}(x)=0\})
\olll(\oLLL'(\tqq)).
\end{align}
Since $X_{h_{1}}$ and $X_{h_{2}}$ are
pairwise independent and $M-0.99$-distributed
on $\R$, from (\ref{*50314a}-\ref{*50314f})
it follows that $\oXX_{h_{1}}$ and $\oXX_{h_{2}}$
are pairwise independent on $\oLLL'(\tqq).$

Put $\oEE_{\delta }=E_{\delta }$.
Clearly, $\oEE_{\delta }$
is periodic by $\ottt$ and
this completes the proof of property (i) in the definition of a $K-M$ family.

It is clear that $\olll(\oEE_{\delta })=
\olll(E_{\delta })<\delta .$
For all $x\not\in \oEE_{\delta }$,
$h\in  \{ 1,...,K  \}$,
let
$\oaaa(x)=\alpha (x)$, $\oooo(x)=
\omega (x)$, $\ottt(x)=\tqq \tau (x)$,
then we have $\frac{\oooo(x)}{\oaaa(x)}>\OOO'\tau (x)=
\OOO \tqq \tau (x)=\OOO \ottt(x)$.

Now we verify \eqref{*12dfw} for $\off_{h}$ and $\oXX_{h}$.
Assume $x\in \oLLL'(\tqq )\sm \oEE_{\delta }$,
$\oaaa(x)=\alpha (x)\leq n<n+m\leq \omega (x)=
\oooo(x),$
$\ottt(x)=\tau (x)\tqq |m$. In fact, it is enough to consider
the case when $\tau (x)\tqq =m=\ottt(x)$.
We claim that for any $h=1,...,K$
we have
$$\frac{1}m
\sum_{k=n}^{n+m-1}
f_{h}(x+k^{2})\leq
\frac{1}m
\sum_{k=n}^{n+m-1}\off_{h}(x+k^{2}),$$
then we will apply \eqref{*12dfw} for $f_{h}$, 
and $X_{h}$.

Since $\tau (x)\tqq =m$
and $x\in \oLLL'(\tqq )\sm \oEE_{\delta }$
implies $x\not\in E_{\delta }$,
using \eqref{**9b} several times we obtain
$$\frac{1}m\sum_{k=n}^{n+m-1}f_{h}(x+k^{2})=
\frac{1}m\sum_{k'=n}^{n+\tau (x)-1}
\sum_{j=0}^{\tqq -1}f_{h}(x+(k'+j\tau (x))^{2})=$$
$$\frac{1}m \sum_{k'=n}^{n+\tau (x)-1}\tqq f_{h}(x+{k'}^{2}).$$
From $x\in\oLLL'(\tqq ),$
it follows that there exists $k_{0}$
such that $x\in -k_{0}^{2}+[i_{k_{0}}\tqq,i_{k_{0}}\tqq +\gamma
2^{\kappa })$, that is,
$x+k_{0}^{2}\in [i_{k_{0}}\tqq
,i_{k_{0}}\tqq + \gamma
2^{\kappa })$ for an $i_{k_{0}}\in\Z,$
and $k_{0}^{2}\in \LLL_{0}'(\tqq )$.
Recall that there are $2^{\kappa }$ many solutions
of $x^{2}\equiv k_{0}^{2}$ modulo $\tqq .$

Since $(\tau (x),\tqq )=1$ for a fixed $k'$, the set
$k'+j\tau (x)$ forms a complete residue system modulo
$\tqq $ as $j$ runs from $0$ to $\tqq -1$, hence there are at least
$2^{\kappa }$ many $j_{k',l}$'s $l=1,...,2^{\kappa }$
with $(k'+j_{k',l}\tau (x))^2\equiv k_{0}^{2}$ modulo $\tqq $,
where $k_{0}$ is defined above. Recalling that
$x\in \oLLL'(\tqq )\sm \oEE_{\delta }$,
for any $k'=n,...,n+\tau (x)-1$, we have
$$\sum_{j=0}^{\tqq -1}
\off_{h}(x+(k'+j\tau (x))^{2})
\geq
\sum_{l=1}^{2^{\kappa }}
\off_{h}(x+(k'+j_{k',l}\tau (x))^{2})
=
$$
$$
\sum_{l=1}^{2^{\kappa }}
f_{h}(x+(k'+j_{k',l}\tau (x))^{2}) \frac{\tqq }{2^{\kappa }}
\geq
 2^{\kappa }f_{h}(x+{k'}^{2})\frac{\tqq }{2^{\kappa }}.$$
Therefore,
$$\frac{1}m\sum_{k=n}^{n+m-1}\off_{h}(x+k^{2})=
\frac{1}m\sum_{k'=n}^{n+\tau (x)-1}\sum_{j=0}^{\tqq -1}
\off_{h}(x+(k'+j\tau (x))^{2})\geq$$
applying \eqref{*12dfw}
for $f_{h}$ and $X_{h}$
$$\frac{1}m\sum_{k'=n}^{n+\tau (x)-1}
\tqq f_{h}(x+{k'}^{2})=
\frac{1}m\sum_{k=n}^{n+m-1}f_{h}(x+k^{2})\geq X_{h}(x)
\geq \oXX_{h}(x).$$
This proves (ii) for $x\in \oLLL'(\tqq)\sm \oEE_{\delta }$.
Since $\oXX_{h}(x)=0$
for
$x\in \R\sm (\oLLL'(\tqq)\cup \oEE_{\delta })$
for these $x$'s \eqref{*12dfw} holds obviously
for $\off_{h}$ and $\oXX_{h}$.

Using \eqref{*14dfww} and $(\tqq,p)=1$ for all
$p\in \ocap$ we have $(\ottt(x),p)=(\tqq\tau (x),p)=1$
and $(\tau \tqq ,p)
=(\ottt,p)
=1$ for all $p\in\ocap .$ This proves
(iii).

To verify (iv), suppose $x\in \oLLL'(\tqq)\sm\oEE_{\delta }$,
$h\in  \{ 1,...,K  \}$ and $\alpha ^{2}(x)=\oaaa^{2}(x)\leq
j<j+\tqq \tau (x)\leq \oooo^{2}(x)=\omega ^{2}(x).$

If $x+j\in \oXXX(\tqq)
$,
then
\begin{align}
\label{*10}
\off_{h}(x+j)&=f_{h}(x+j)\tqq /2^{\kappa }=
f_{h}(x+j+\tau(x))\tqq /2^{\kappa }=...=\\
&f_{h}(x+j+\tqq \tau(x))\tqq /2^{\kappa }=
\nonumber\off_{h}(x+j+\tqq \tau (x))
\end{align}
when $\alpha ^{2}(x)=\oaaa^{2}(x)
\leq j < j+\tqq \tau (x)\leq \oooo^{2}(x)=
\omega ^{2}(x)$.

If $x+j\not \in \oXXX(\tqq)
$, then
$\off_{h}(x+j)=0=\off_{h}(x+j+\tqq \tau(x)).$
This verifies (iv).

Next we prove (v)
\begin{equation}\label{*11b}\frac1{\tau \tqq }
\int_{0}^{\tau \tqq }\off_{h}=
\ointt \off_{h}\leq \GGG\cdot \gamma \cdot  2^{-M+1}.
\end{equation}
Indeed,
$$\frac1{\tau \tqq }
\sum_{j=0}^{\tqq -1}
\int_{0}^{\tau }\off_{h}(x+j\tau )dx=
\frac1{\tau }\int_{0}^{\tau }
\frac{1}{\tqq }
\sum_{j=0}^{\tqq -1}\off_{h}(x+j\tau )dx=(*).$$
To continue this computation recall
that
from $(\tau ,\tqq )=1$ it follows that
 $\lf x \rf+j\tau $ hits each residue class modulo
$\tqq $ once
as $j$ varies from $0$ to $\tqq -1$
and $f_{h}(x+j\tau )=f_{h}(x)$
for all $j$. Thus, recalling that
$\gamma '$ associated to $\Lambda=\R$ equals $1$
and using (v) for $f_{h}$
$$(*)=
\frac{1}\tau \int_{0}^{\tau }\frac1{\tqq }
f_{h}(x)\frac{\tqq }{2^{\kappa }}\gamma 2^{\kappa
}dx=
\frac{1}{\tau }\gamma \int_{0}^{\tau }
f_{h}(x)dx\leq \GGG\gamma 2^{-M+1}.$$
This proves \eqref{*11b}.
\end{proof}

\subsection{Proof of Lemma \ref{lksystem}}\label{42}

Let a positive integer $M$ be given together with input parameters
$\delta > 0, \Omega, \Gamma < 1, A \in \N$ and $\capp\sse \N$ with
sufficiently large complement. To show the existence
of a $(K+1)-M$ family living on $\R$, we need to fix several constants for
the induction argument.

To begin with we will use the following Lemma  \ref{rdmlem} to
choose a ``leakage constant''
which will remain
fixed during the inductive construction of a $(K+1)-M$ family
from $K-M$
families.
This lemma is a more exact expression of the ideas given in
Remark \ref{**5dfw}.

Considering
the sets $\oLLL'_{\gamma }(q)$, a direct calculation shows: 
\begin{equation}\label{*19bud}
\olll(\oLLL'_{\gamma }(q))<\gamma
\text{ and }\olll(\R\sm \oLLL'_{\gamma }(q))>1-\gamma .
\end{equation}
 However,
the closer $\gamma $ to $0$, the smaller the percentage
of ``loss due to overlaps". To obtain estimates from the opposite sides
we will use Lemma \ref{rdmlem}.

\begin{lemma}\label{rdmlem} For each $0<\gamma <1/7$, one can choose constants
$C_{\gamma }>1>\tCC_{\gamma }>0$, $\kappa _{\gamma }\in \N,$ such that
for each $\kappa >\kappa _{\gamma }$ there exists
$p_{\gamma,\kappa }$ for which  if $p_{\gamma,\kappa }<p_{1}<...
<p_{\kappa }$ and  $q=p_{1}\cdot \cdot \cdot p_{\kappa }$,
then
\begin{equation}\label{*13/1}
C_{\gamma }>\frac{\gamma }{\olll(\oLLL'_{\gamma }(q))} \ \ \
\text{ and } \ \ \ \
\tCC_{\gamma }\olll(\R\sm \oLLL'_{\gamma }(q))<1-\gamma -\gamma ^{2}.
\end{equation}
In fact, we can choose
$$ C_{\gamma } =
\frac{1} {1-7\gamma}\text{ and }\tCC_{\gamma } = \frac{1-\gamma -
 \gamma^2}{1-\gamma +7\gamma^2}.$$
Therefore, $C_{\gamma }\to 1$ and $\tCC_{\gamma }
\to 1$ as $\gamma \to 0+.$
\end{lemma}

We remark that in \eqref{*13/1} the second order term
in $(1-\gamma -\gamma ^{2})$ appears for technical reasons.
It is clear that $(1-\gamma -\gamma ^{2})/(1-\gamma )\to 1$
as $\gamma \to 0+.$

\begin{proof}
We use the fact that the limiting
 distribution of the gaps between squares is continuous. In fact, consider
$q = p_1\cdot...\cdot p_\kappa$ where $p_1 < ... < p_\kappa$ are
odd primes. Let $0 = x_1 < x_2< ...< x_{\sigma_q} < q  = x_{\sigma_q +1}$
be the squares mod $q$
so that $\sigma_q = \prod_1^\kappa (\frac{p_i+1}{2})$.
Let each gap
 $g_i = x_{i+1}-x_i$
have weight $1/\sigma_q$.
The expected gap size is
$s_q = \frac{1}{\sigma_q}\sum_{1}^{\sigma_q}g_i =\frac{q}{\sigma_q} =
2^\kappa\prod_1^\kappa \ (\frac{p_i}{1+p_i})$.
Let us normalize the gaps; $\ y_i
 =\frac{g_i}{s_q}$.
Kurlberg and Rudnick in \cite{[Kur]}, Lemma 14
proved the following result. For each $x \in \R$,
\begin{equation}\label{*23zz}
\lim_{\kappa\to \infty}
\frac{\#\{i: \ y_i  \ \leq \ x\}}{\sigma_q} = 1 -e^{-x}.
\end{equation}
Choose $\kappa _{\gamma}$ such that if
$\kappa >\kappa _{\gamma }$, then
$\frac{\#\{i: \ y_i  \ \leq \ 2\gamma \}}{\sigma_q} < 1-e^{-3\gamma }
<3\gamma $.
For each $\kappa >\kappa _{\gamma }$ choose
$p_{\gamma,\kappa }$ such that  if $p_{\gamma,\kappa }<p_{1}<...
<p_{\kappa}$ then $\prod_1^\kappa (1+\frac{1}{p_i}) <2$ and
$\prod_1^\kappa (1-\frac{1}{p_i}) > 1-\gamma.$ Let
$q=p_{1}\cdot \cdot \cdot p_{\kappa }$.
Letting $x_i'$ be the squares modulo $q$ in $[0,q)$ which are not
divisible by any of the prime factors of $q$, we have
\begin{align}
\olll(\oLLL'_{\gamma }(q)) &\geq \frac{\gamma2^\kappa}{q}\#\{j:x_{j+1}'-x_j'
> \gamma2^\kappa\} \\
&\geq \frac{\gamma2^\kappa}{q}
\left(\prod_1^\kappa (\frac{p_i -1}{2})-\#\{i:y_i \leq
\frac{\gamma2^\kappa}{s_q} = \gamma \prod_1^\kappa (1+\frac{1}{p_i}) < 2
\gamma\}\right).
\end{align}
By \eqref{*23zz} and our assumptions, we get
\begin{align}\label{*080211}
\olll(\oLLL'_{\gamma }(q)) &\geq \gamma\left((1-\gamma)
-{3\gamma }\prod_1^\kappa(1+\frac{1}{p_i})\right )> \\
\nonumber & \gamma((1-\gamma) -3\gamma\cdot 2 ) = \gamma(1-7\gamma)
= \frac{\gamma}{C_\gamma}.
\end{align}
We have
$$
\olll(\oLLL'_{\gamma }(q)) > \gamma(1-7\gamma) = 1 -
\frac{1-\gamma-\gamma^2}{\tCC_{\gamma }}.
$$
So,
$$
\frac{1-\gamma-\gamma^2}{\tCC_{\gamma }} > 1 - \olll(\oLLL'_{\gamma
  }(q)) = \olll(\R\sm \oLLL'_{\gamma }(q)).
$$

\end{proof}

In order to apply
Lemma \ref{lputres} we
need to choose a
positive
``leakage constant," $\gamma $,
which remains fixed during all steps of the leakage
producing the $(K+1)-M$ family.

{\it Fixing the leakage constant} $\gamma.$\\
 We
choose
$0<\gamma _{0}<10^{-7}$ so that
\begin{equation}\label{*24zz}
C_{\gamma _{0}}=\frac{1}{1-7\gamma _{0}}<\Gamma.
\end{equation}
Moreover,
for each $\gamma <\gamma _{0}$ with
$\gamma =2^{-c_{\gamma }}$ where $c_{\gamma }\in \N$,
by Lemma \ref{rdmlem}
we choose $\kappa _{\gamma }$ such that
for all $\kappa >\kappa _{\gamma }$
there exists $p_{\gamma,\kappa  }$ for which if
$p_{\gamma,\kappa  }
<p_{1}<...<p_{\kappa }$, $q=p_{1}\cdot \cdot \cdot p_{\kappa },$
then
\begin{equation}\label{*15d}
\olll(\oLLL'_{\gamma }(q))>
\frac{\gamma }{C_{\gamma _{0}}}=
(1-7\gamma _{0})\gamma
>
\frac{9\gamma }{10}\text{ and }
\olll(\R\sm \oLLL'_{\gamma }(q))<1-\frac{9\gamma }{10}.
\end{equation}
We also have,
$$\tCC_{\gamma }>1-\gamma _{0}-\gamma _{0}^{2}>1-10^{-6}.$$

From now on a value of $\gamma <\gamma _{0}$,
with $\gamma =2^{-c_{\gamma }}$, $c_{\gamma }\in \N$
satisfying the above assumptions is fixed.

We note that the only input parameter that the leakage constant
depends on is $\Gamma.$ 

We will write $\LLL(q)$, $\oLLL(q)$ and $\oLLL'(q)$
instead of
$\LLL_{\gamma }(q)$, $\oLLL_{\gamma }(q)$ and $\oLLL'_{\gamma }(q)$,
respectively.

Next,
after giving an outline
we start the details of the proof of
 Lemma \ref{lksystem}.

\subsubsection{Setting up the induction argument for Lemma
\ref{lksystem}}\label{421}
\begin{proof}
We proceed by mathematical induction.

To start our induction we need to show that $1-M$ families
exist. During the general step of our induction we show
that from the existence of $K-M$ families one can deduce
the existence of $(K+1)-M$ families.
Since many steps of the $1-M$ family case are
shared with the general $K-M$ family case we work
out our argument so that it can be used for the
later induction steps
without any unnecessary duplication.
Therefore, working on the first step of our induction one
should think of $K=0$ during the first reading of Sections
\ref{422}-\ref{428} and obtain this way the $(K+1)-M$, that
is, the $1-M$ families. Then in Sections \ref{422b}-\ref{428b}
we discuss the alterations needed for $K>0$.
It will be useful to keep in mind that if $K=0$ then
$f_{K+1,0}=f_{1,0}$, only $h=1\in  \{ 1,...,K+1  \}$
and there is no $h\in  \{ 1,...,K  \}$.

Now we discuss briefly our general plan.
When $K>0$
we assume that $K-M$ families living
on $\R$ exist
for all possible input parameter choices.
Let
$\delta >0,$ $\Omega,$ $\GGG>1$, $A\in \N$
and $\capp\sse \N$ be given.
We will define our $(K+1)-M$ family with
these
input parameters.
We can assume that $\capp$ is closed under products.
During the definition
of the $(K+1)-M$ family another, ``inner" finite induction is used
(with respect to $L$)
which is called the {\it leakage process}.
This technically delicate process is the focus
of the next several sections.
During this process we will
use Lemma \ref{lputres} to
define families which are almost
$(K+1)-M$
families on sets of the type $\oLLL'(q)$,
except for the new functions $f_{K+1,L}$
and $X_{K+1,L}$.
Each $f_{K+1,L}$ is the indicator function
of a set.
As $L$ grows the support of $f_{K+1,L}$
decreases.
Lemma \ref{probabilitylemma}
and Lemma \ref{taurearrangement}
are used to ensure that ``squares hit sufficiently often"
the support of $f_{K+1,L}$. (This motivated the term
``leakage" since the values of $f_{K+1,L}$ leak onto
some larger sets when we consider averages along the squares.
See also \cite{BM}.)
This also requires that before defining $f_{K+1,L}$
one uses Lemma \ref{taurearrangement} to choose
$\tau _{L-1}'$ and make a $\tau _{L-1}'$
rearrangement to yield an intermediate
function $f_{K+1,L-1}'$. At the same time we must keep
track of our new random variable $X_{K+1,L}$ and
other auxiliary functions. It is essential in this
induction that we can vary $\tau$ and $\kappa.$

To help the reader going through the details of
the proof here is an outline of the main features
of the various sections of the proof. We hope
this outline might help prevent the reader
from becoming lost in the details of the proof.

When $K=0$ in Section \ref{422} we start
the leakage process.
We define $f_{K+1,0}\equiv 1$ and $X_{K+1,0}$.
At this stage $f_{K+1,0}$ is supported on $\R$.
During the leakage process the size of the support
of the functions $f_{K+1,L}$ is shrinking and we
are interested in how much of $f_{K+1,L}$
is ``leaking" onto larger sets.
When $K>0$ in Section \ref{422b},
in addition, we introduce
a $K-M$ family
periodic by $\tau_{0} $, consisting
of functions $f_{h,0}$, $X_{h,0}$ for
$h=1,...,K$.

In Section \ref{423} (see also Section \ref{423b} when $K>0$)
we assume that we have accomplished
step $L-1$ of the leakage and we have a  family
periodic by $\tau _{L-1}$, consisting of functions
$f_{h,L-1}, $ $X_{h,L-1}$ for $h=1,...,K+1$.
We also introduce the auxiliary sets $S_{L-1,l}$,
$l=0,...,L-1$ used to describe the distribution
of $X_{K+1,L-1}$.

In Section \ref{424} (see also Section \ref{424b} when $K>0$)
we choose a prime number
$\tau _{L-1}'$ which is much larger than $\tau _{L-1}$
and by using Lemma \ref{taurearrangement} we perform
a $\tau '_{L-1}$ rearrangement of the family
coming from Section \ref{423}. This way we obtain
a family periodic by $\tau '_{L-1}$, consisting of functions
$f_{h,L-1}'$, $X_{h,L-1}'$, $h=1,...,K+1$.
The auxiliary sets
used for describing the distribution of
$X_{K+1,L-1}'$ are denoted by $S_{L-1,l}'$, $l=0,...,L-1.$

In Section \ref{425}
we choose a $\kappa _{L}$ and
a square free number $q_{L}=p_{1,L}\cdot \cdot \cdot q_{\kappa _{L},L}$
such that $2^{\kappa _{L}}$ is much larger than
$\tau '_{L-1}$. The average value of the difference
between elements of $\LLL_{0}(q_{L})$ is close to
$2^{\kappa _{L}}$. We introduce
some auxiliary sets, among them $\FFF_{L}$ and
$\PPP_{L}$, so that $\FFF_{L}\sse \R\sm \oLLL'(q_{L})
\sse \PPP_{L}$ and these two auxiliary sets
consist of intervals of the form
$[j\tau '_{L-1},(j+1)\tau '_{L-1}),$ $j\in \Z.$
If $2^{\kappa _{L}}$ is much larger than
$\tau _{L-1}'$, then $\olll(\FFF_{L})$ and
$\olll(\PPP_{L})$ both approximately equal
$\olll(\R\sm \oLLL'(q_{L})).$
To define our $(K+1)-M$ family on $\PPP_{L}$ (which
is approximately $\R\sm \oLLL'(q_{L})$)
we will use mainly the functions coming from Section
\ref{424}.
We define $X_{h,L}(x)=X_{h,L-1}'(x)$
if $h\leq K+1$ and $x\in \FFF_{L}$.
This Section is identical for the cases $K=0$
and $K>0$.

For $K>0$ in Section \ref{426b} by using Lemma \ref{lputres} we put
a $K-M$ family onto $\oLLL'(q_{L})$.
This will
yield functions $\off_{h,L}$, $\oXX_{h,L}$
periodic by $\ottt_{L} q_{L}$.
For $h=1,...,K$ we define $X_{h,L}$ on
$\oLLL'(q_{L})$ by using $\oXX_{h,L}$.
For $h=1,...,K$ our functions will be sums
of $f'_{h,L-1}$ restricted to $\FFF_{L}$ 
and of
the functions
$\off_{h,L}$ ``living" on $\oLLL'(q_{L}).$
This combined family will be
 periodic by $\tau _{L}=q_{L}\ottt_{L}\tau _{L-1}'.$
The $K=0$ version discussed in Section
\ref{426} is much simpler  because we do not have to deal
with this putting a lower level family on  the 
$\oLLL'(q_{L})$ step.

The ``leakage" is done when we define $f_{K+1,L}$
so that it equals the restriction of
$f'_{K+1,L-1}$ onto the set $\PPP_{L}$.
This means that the support $F_{L}$
of $f_{K+1,L}$ will have a very small intersection
with $\oLLL'(q_{L})$. ``Most" of $F_L$
will be a subset of $\R\sm \oLLL'(q_{L})$
and will approximately equal the auxiliary
set $S_{L,0}$. The nested sequence
$S_{L,l}$, $l=0,...,L$ will describe
the distribution $X_{K+1,L}$, the larger
$l$, the smaller the values
 $X_{K+1,L}$
can take
on $S_{L,l}\sm S_{L,l-1}$.

In Section \ref{427} we make the calculations needed
to show that we have enough ``leakage"
from the support of $f_{K+1,L}$ so that
we have
the domination inequality
\eqref{*12dfw} with $f_{K+1,L}$ and
$X_{K+1,L}$. This section is again the same
for the cases $K=0$ and $K>0$.

Finally, in Section \ref{428}
(see also Section \ref{428b} when $K>0$)
 we terminate the leakage
process when we have reached a suitably large
$L=L''\leq L'$.
The functions $f_{h,L''}$ for $h=1,...,K+1$ will
yield the functions
$f_{h}$
we need for the $(K+1)-M$ family.
For $h=1,...,K$ the functions
$X_{h}$ of the $(K+1)-M$ family will equal the functions
$X_{h,L''}$. To define $X_{K+1}$ we use the sets $S_{L'',l}$,
$l=0,...,L''$ related to the distribution of
$X_{K+1,L''}$. We will choose $X_{K+1}$ so that it is
$M$$-$$0.99$ distributed and less or equal than $X_{K+1,L''}$.

\medskip

Before turning to the details of the induction
to help the reader going through the details of the proof
for easy reference we collect some defintions and properties (some
of them will be discussed later during the proof) at the same place.\\
{\Large {\bf Quick reference  summary:}}\\
$\LLL(q)=\LLLg(q)=-\LLL_{0}(q)+ \{ j\in\Z:
0\leq j<\gamma 2^{\kappa }  \},$\\
$\oLLL(q)=\oLLLg(q)=\LLLg(q)+[0,1)=
-\LLL_{0}(q)+ \{ x:0\leq x< \gamma 2^{\kappa }  \},$\\
$\LLL_{0}'(q)= \{ n\in \LLL_{0}(q): p_{j}\not|n, \text{ for all
}j=1,...,\kappa  \},$\\
$\LLL'(q)=\LLL'_{\gamma  }(q)
=-\LLL'_{0}(q)+ \{ j\in\Z:0\leq j <\gamma 2^{\kappa }  \},$\\
$\oLLL'(q)=\oLLL_{\gamma  }'(q)=\LLL_{\gamma  }'(q)+[0,1)=-\LLL_{0}'(q)+
 \{ x:0\leq x <\gamma 2^{\kappa }  \}.$
 
{\bf At Step $L=0$ of the leakage process we have:}
$f_{K+1,0}\equiv 1$,
$F_{0}=\R,$
$S_{0,0}=\R$, $r_{0}=1,$
$X_{K+1,0}\defeq
(1-\rho ')\tCC_{\gamma }=
(1-\rho ')\tCC_{\gamma }\olll(F_{0})<1$.

{\bf After Step $L-1$ of the leakage process we have:}
The set $F_{L-1}$,
an exceptional set
$E^{L-1}$, a period $\tau_{L-1}$
such that $F^{L-1}$,
$E^{L-1}$,
$X_{h,L-1},$ $f_{h,L-1}$, $(h=1,...,K+1)$ are
periodic by $\tau_{L-1} ,$
$f_{h,L-1}:\R\to [0,\oo)$,
 the ``random" variables
 $X_{h,L-1}:\R\to\R$
are pairwise independent  for $h=1,...,K+1,$
$X_{h,L-1}$ are $M$$-$$0.99$-distributed
 for $h=1,...,K$, but not for $h=K+1.$
For the distribution of $X_{K+1, L-1}$ 
the auxiliary sets $S_{L-1,l}$ are used.
For all $x\not\in E^{L-1}$
there exist
$\omega_{L-1} (x)>\alpha_{L-1} (x)>A,$ $\tau_{L-1} (x)< \tau_{L-1}.$

{\bf Then we do a $\tttt'_{L-1}$ periodic rearrangement.}
We choose and fix a sufficiently large
prime
$\tau '_{L-1}$. The set 
$F_{L-1}^{\tau _{L-1}'}$ is
the $\tau _{L-1}'$ periodic rearrangement of
$F_{L-1}$.
We modify our sets and functions so that they are all
periodic with respect to $\tau _{L-1}'$.
The new functions are:
$f'_{h,L-1}$
$\alpha '_{L-1}(x)$
$\omega _{L-1}'(x)$
$\tau _{L-1}'(x)$
$X_{h,L-1}'(x)$.
We define ${E'}^{L-1}$ so that it satisfies \eqref{*h16}.
We set $S'_{L-1,L-1}=\R$ and
 define the sets $S'_{L-1,l}$.
 These sets are used in \eqref{*16b} for the distribution of
 $X_{K+1,L-1}'$.

{\bf Then we choose a sufficiently large $q_{L}$,
and several important auxiliary sets:}
$\oXXX(q_{L})=
\cup_{j\in \Z} [jq_{L} ,jq_{L} +\gamma 2^{\kappa_{L} })$,\\
$\tFFF_{L}= \{ x:dist(x,\oLLL'(q_{L})\cup
\oXXX(q_{L}))>2\tau '_{L-1}  \},$\\
$
\FFF_{L}=
\cup \{ [j\tau '_{L-1},
(j+1)\tau '_{L-1}):
\tFFF_{L}\cap
[j\tau '_{L-1},(j+1)\tau '_{L-1})\not=\ess
\},
$
which satisfy:
$\tFFF_{L}\sse\FFF_{L}\sse \hFFF_{L}\sse \R\sm \oLLL'(q_{L})
\sse \hPPP_{L}\sse \PPP_{L}
\sse \tPPP_{L}$
The sets
$\FFF_{L}$ and
$\PPP_{L}$ are periodic by $\tau '_{L-1} q_{L}$
and the sets $\tFFF_{L}$, $\hFFF_{L}$, $\tPPP_{L}$,
and $\hPPP_{L}$ are periodic by $q_{L}$.
They satisfy \eqref{*16}.
We put 
$E_{L}'=\PPP_{L}\cap {E'}^{L-1}$,
$E_{L}''=\R\sm (\oLLL'(q_{L})\cup \tFFF_{L})$
and $\tEE_{L}''=
\R\sm (\oLLL'(q_{L})\cup \FFF_{L})\sse E_{L}''$.
The set
 $E_{L}''$ is periodic by $q_{L}$, while $\tEE_{L}''$ is periodic
by $\tau _{L-1}' q_{L}$.
The exceptional set $\oEE_{L}'''$ is defined at \eqref{*EvvL}
and $E'''_{L}$ in \eqref{*EvvL2}.

{\bf From this place on our definition of our new objects, like
the functions
$f_{h,L}$ and $X_{h,L}$ splits and follows two different paths.
One will be the definition of these objects on
$\R\sm \oLLL'(q_{L})$ and the other the definition
of these objects on $\oLLL'(q_{L})$.}

For the first path we are unable to use exactly
the set $\R\sm \oLLL'(q_{L})$. We need to use 
the auxiliary sets $\FFF_{L}$ and $\PPP_{L}$ which are
good approximations of this set. We have for example
$f_{K+1,L}=f'_{K+1,L-1}\chi_{\PPP_{L}}$.

For the second path by our induction assumption we put a 
$K-M$-family on $\oLLL'(q_{L})$. The sets and functions obtained at this step
are periodic by $\ottt_{L}q_{L}.$ They are $\off_{h,L},$ $\oaaa_{L},$
$\oooo_{L},$ $\ottt_{L}$, $\oXX_{h,L}$ and there is an exceptional set $E_{\ddd,L}.$

We combine these two paths 
when we define
$f_{h,L}=f'_{h,L-1}\cdot\chi _{\FFF_{L}}
+\off_{h,L} \text{ for }h=1,...,K.$

{\Large {\bf Quick reference summary ends here.}}

\medskip

Next we turn to the details of our argument.
We start by choosing some constants.

Choose a positive integer $L'$ such that
\begin{equation}\label{*h6}
(1-\frac{\gamma }2)^{L'}<\frac{1}{2^{M}}.
\end{equation}

By recursion we construct pairs of functions
$f_{K+1,0},X_{K+1,0},...,f_{K+1,L},X_{K+1,L}$ and some associated objects.
The inner finite induction, the ``leakage"
will halt at some step $L''\leq L'.$ We show that the functions
$f_{K+1,L''}, X_{K+1,L''}$ and their associated objects form a $(K+1)-M$ family
except for the distribution of $X_{K+1,L''}.$ However we will know enough
about its distribution to easily obtain a $(K+1)-M$ family.

By \eqref{*24zz} we have $C_{\gamma }<C_{\gamma _{0}}<\GGG$.

We recall that the only input parameter that the leakage constant
depends on is $\Gamma.$ This dependence, and the possibility
of using different values for $\Gamma$ will play an important
role during the definition of $(K+1)-M$ families with $K>0$.

Fix a constant
$\GGG_{0}>1$ such that
\begin{equation}\label{*16h}C_{\gamma }\GGG_{0}<\GGG.
\end{equation}
Put
$$\delta _{L}=\frac{\delta }{4(L'+1)}\text{ for }
L=0,...,L'.$$

Next we  choose
 sufficiently small positive constants
$\rho $, $\rho '$ and $\trrr$.

We suppose that \begin{equation}\label{*39x}
(1-\rho )(1-2\gamma )>1-3\gamma.
\end{equation}

Recall that $1>\tCC_{\gamma}> 1-10^{-6}$ and we choose
$\rho '>0$ such that
\begin{equation}\label{*30x}
0.999<\tCC_{\gamma }(1-\rho ')<1<1.001.
\end{equation}

Since $0<\gamma <\gamma _{0}< 10^{-7}$ and
$1>\tCC_{\gamma}>1-10^{-6}$ we can suppose that
$\rho$ and $\rho'$ are so small that
\begin{equation}\label{*108x}
(1-\rho )^{2}(1-2\gamma )^{2}(1-\rho ')\tCC_{\gamma }-
\frac{(1-\rho ')\tCC_{\gamma }}{0.999\cdot (1-\rho )^{2}}
\cdot \frac{1}2>0.99\cdot \frac{1}2.
\end{equation}

Moreover, choose $\trrr>0$ such that
\begin{equation}\label{*38x}
(1-\frac{\rho '}2)(1-2\trrr)>1-\rho '.
\end{equation}

Finally, we set $\capp_{0}=\capp
\cup\capp_{0}'\cup...\cup\capp_{L'}'$
where $\capp$ and each $\capp_{j}'$ contains infinitely
many primes and all their possible products, but numbers
in different sets are relatively prime; moreover
$\capp_{0}$ has sufficiently large complement.

\subsubsection{Step $L=0$ of the leakage process}\label{422}

We put $f_{K+1,0}\equiv 1$. Set
$F_{0}=\R,$
$S_{0,0}=\R$, $r_{0}=1.$
So, $\olll(F_{0})=1$ and  $X_{K+1,0}\defeq
(1-\rho ')\tCC_{\gamma }=
(1-\rho ')\tCC_{\gamma }\olll(F_{0})<1$, see \eqref{*30x}.

For the case $K=0$ we use the following argument.
For $K>0$ see a different argument in \ref{422b}.

We
choose a sufficiently large $\tau_{0}$, and functions
$\alpha _{0}(x)$, $\omega _{0}(x)$, $\tau_{0} (x)$ taking
integer values for all $x\in\R$ (in fact, these functions can be
constant on $\R$), such that the following assumptions hold:
$\omega_{0} (x)>\alpha_{0} (x)>A,$ $\tau_{0} (x)< \tau_{0} $,
$\omega_{0} ^{2}(x)<\tau_{0} ,$
$\frac{\omega_{0} (x)}{\alpha_{0} (x)}>\OOO \tau_{0} (x)
=\OOO_{0} \tau_{0} (x)$,
and
for all $p\in\capp_{0}$, $(\tau_{0} (x),p)=1,$ $(\tau_{0} ,p)=1.$
For example, we could take $\alpha _{0}(x)=A+1$, $\tau _{0}(x)$
to be the smallest odd prime which is relatively prime to
all elements of $\capp_{0}$, $\omega _{0}(x)=\tau _{0}(x)
\OOO (A+2)$ and $\tau _{0}$ be the smallest prime relatively
prime to the elements of $\capp_{0}$ and greater than $(\omega_{0}(x))^2
 $.

By our choices; $f_{K+1,0}=f_{1,0}\equiv 1$.
 We put $E_{\delta _{0}}=\ess$
and $E^{0}=\ess.$
For any $m\in \N$  we have the all important ``domination'' property:
for all $x \in \R\setminus E^0$,
$$
\frac{1}m\sum_{k=n}^{n+m-1}f_{1,0}(x+k^{2})>X_{1,0}(x).$$
It is also clear that
for all $x\in\R$,
$f_{1,0}(x+j+\tau_{0} (x))=f_{1,0}(x+j)$
for any $j\in \R$.

\subsubsection{The setting after step $L-1$ of the leakage}\label{423}

Assume we have accomplished step $L-1$ of the leakage
process.
We have constructed some objects satisfying the
following conditions.
There is an exceptional set
$E^{L-1}$
with
\begin{equation}\label{*h10}
\olll(E^{L-1})<\frac{L}{(L'+1)}\delta;
\end{equation}
$\capp_{L-1}=\capp\cup\capp_{L}'\cup...\cup \capp_{L'}';$
there exists a period $\tau_{L-1}$
such that
$E^{L-1}$,
$X_{h,L-1},$ $f_{h,L-1}$, $(h=1,...,K+1)$ are
periodic by $\tau_{L-1} ,$
$f_{h,L-1}:\R\to [0,\oo)$,
(for $K>0$) the ``random" variables
 $X_{h,L-1}:\R\to\R$
are pairwise independent  for $h=1,...,K+1,$
$X_{h,L-1}$ are $M$$-$$0.99$-distributed
 for $h=1,...,K$ (in (\ref{**h12}-\ref{*15b})
we list the assumptions
about the distribution of
$X_{K+1,L-1}$, recall that for  $K=0$ there is no $h$ satisfying
$h=1,...,K$).
For all $x\not\in E^{L-1}$
there exist
$\omega_{L-1} (x)>\alpha_{L-1} (x)>A,$ $\tau_{L-1} (x)< \tau_{L-1} $
such that
$\omega_{L-1} ^{2}(x)<\tau_{L-1} ,$
$\frac{\omega_{L-1} (x)}{\alpha_{L-1} (x)}>\OOO \tau_{L-1} (x)$;
moreover
if $\alpha_{L-1} (x)\leq n<n+m\leq \omega_{L-1} (x)$
and $\tau_{L-1} (x)|m$, then for all $h=1,...,K+1$,
\begin{equation}\label{*c17}
\frac{1}m\sum_{k=n}^{n+m-1}f_{h,L-1}(x+k^{2})>X_{h,L-1}(x);
\end{equation}
for all $p\in\capp_{L-1}$, $(\tau_{L-1} (x),p)=1,$ $(\tau_{L-1} ,p)=1;$
moreover for all $x\not\in E^{L-1}$,
$f_{h,L-1}(x+j+\tau_{L-1} (x))=f_{h,L-1}(x+j)$ whenever
$\alpha_{L-1} ^{2}(x)\leq j<j+\tau_{L-1} (x)\leq \omega_{L-1} ^{2}(x)$
for all $h\in  \{ 1,...,K+1  \}$.

We suppose
that the values of $f_{K+1,L-1}$ are $0$ or $1$, that is,
it is an indicator function.

If $L-1=0$, then $X_{K+1,L-1}$ is constant.

If $L-1>0$,
that is, $L\geq 2$
then we give the extra assumptions about the distribution of
$X_{K+1,L-1}$ as follows.

Recall $F_{0}=\R$ and also recall
that for a Lebesgue measurable set $F$, periodic by
$p$ we have $\olll(F)=\frac{1}p\lambda (F\cap [0,p))=
\lim_{N\to \oo}\frac{\lambda(F\cap [-N,N])}{2N}.$
We suppose that the sets $F_{l}$ periodic by
$\tau_l$ and the numbers
$r_{l}$ have been defined for $l=0,...,L-2$ during
the previous steps of our induction,
\begin{equation}\label{*www9}
\frac{\olll(F_{l})}{\olll(F_{l-1})}=r_{l}\text{ and }1-2\gamma <r_{l}
<1-\frac{\gamma }2
\end{equation}
hold for $l=1,...,L-1.$ Clearly, $\olll(F_l)=r_0\cdots r_l.$
We also have $\tau_l$ periodic functions $f_{K+1,l} = \chi_{F_l}$
for each $0 \leq l \leq L-1$. 

Set $F_{L-1}= \{ x:f_{K+1,L-1}(x)=1  \}$ and
$r_{L-1}=\frac{\olll(F_{L-1})}
{\olll(F_{L-2})}<1.$ We also assume
\begin{equation}\label{*h12}
1-2\gamma <r_{L-1}<1-\frac{\gamma }2,
\end{equation}
$\olll(F_{L-1})=r_{1}\cdot \cdot \cdot r_{L-1}
=r_{0}\cdot \cdot \cdot r_{L-1}.$
 In \eqref{50401} we explicitly show that this holds
for $r_1$. The sets $S_{L-1,0}\sse ...\sse S_{L-1,L-1}=\R$
are defined so that
\begin{equation}\label{**h12}
\frac{1}{1-\rho}
\olll(F_{L-1})>\olll(S_{L-1,0})>(1-\rho )\olll(F_{L-1}),
\end{equation}
if $x\in S_{L-1,0}$ then
$X_{K+1,L-1}(x)=(1-\rho ')\tCC_{\gamma }=(1-\rho ')
\tCC_{\gamma }\olll(F_{0}).$
For $l=0,...,L-1$
we have
\begin{equation}\label{*h13}
\frac1{(1-\rho )r_{0}\cdot \cdot \cdot r_{l}}
\olll(S_{L-1,0})>\olll(S_{L-1,l})>
\frac{1-\rho }
{r_{0}\cdot \cdot \cdot r_{l}}\olll(S_{L-1,0}),
\end{equation}
which is equivalent to
\begin{equation}\label{**h13}
\frac1{(1-\rho )} \olll(S_{L-1,0})>
\olll(F_{l})\cdot \olll(S_{L-1,l})>
(1-\rho) \olll(S_{L-1,0}).
\end{equation}
If $x\in S_{L-1,l}\sm S_{L-1,l-1}$ for $l\in  \{ 1,...,L-1  \},$
then
\begin{equation}\label{*15b}
X_{K+1,L-1}(x)=(1-\rho' )r_{0}\cdot \cdot \cdot r_{l}
\tCC_{\gamma }=(1-\rho ')\olll(F_{l})\tCC_{\gamma }.
\end{equation}
The sets $S_{L-1,l}$ are increasing almost by a factor
$1/r_{l}$ in size, whereas the value of
$X_{K+1,L-1}$ on the difference is decreasing
by a factor $r_{l}$.
We also assume that $F_{L-1}$ has the property that if
$x\in F_{L-1}$ then $[\lf x \rf, \lf x \rf+1)\sse F_{L-1}.$

We note that by \eqref{*www9}
$$
\ointt f_{K+1,l} \leq (1 -\frac{\gamma}{2})^l\text{ for }l=0,...,L-1.
$$

\subsubsection{Rearrangement with respect to $\tau '_{L-1}$,
choice of $\tau '_{L-1}$}
\label{424}

In order to construct the next set of objects in the recursion,
we first create, by
 rearrangement, some associated objects to the  $(L-1)$st  step
which are denoted by
attaching primes.

Since the set $F_{L-1}$ is periodic by $\tau _{L-1}$
and is the union of some integral intervals we can
apply Lemma \ref{taurearrangement}.
We choose $M_{\rho'/2}$ such that for all
prime numbers
$\tau '_{L-1}>
M_{\rho '/2}$ if we consider
$F_{L-1}^{\tau _{L-1}'}$,
the $\tau _{L-1}'$ periodic rearrangement of
$F_{L-1}$,
then for any $x\in \R$ if $\tau _{L-1}'|m$,
then
\begin{equation}\label{*14}
\frac{1}m\sum_{k=n}^{n+m-1}\chi _{F_{L-1}^{\tau _{L-1}'}}
(x+k^{2})\geq (1-\frac{\rho '}2)\olll(F_{L-1}).
\end{equation}
We will choose and fix a sufficiently large
prime
$\tau '_{L-1}\in \capp_{L-1}'$. We define the numbers
$\cat_{i}$, $i=1,...,5$ below.
We choose $\tau '_{L-1}$ so that
it is larger than the maximum of
$\cat_{1}$,...,$\cat_{5}$ and
$M_{\rho '/2}$.
Hence 
\eqref{*14},
\eqref{*29x},
\eqref{*h17}, \eqref{**h17}, \eqref{***h17},
\eqref{*85xx}, \eqref{*19d}, and \eqref{*h18}
hold.

Now, we modify our sets and functions so that they are all
periodic with respect to $\tau _{L-1}'$.
Since we are going to define functions which are periodic by
$\tau _{L-1}',$ it is sufficient to define them
on $[0,\tau _{L-1}').$

If $x\in[0,\lf \tau '_{L-1}/\tau_{L-1} \rf \cdot \tau _{L-1})$
and the right hand side of the equation is defined at $x$,
set\\
\begin{align*}
f'_{h,L-1}(x)&=f_{h,L-1}(x),\ h=1,...,K+1,\\
\alpha '_{L-1}(x)&=\alpha _{L-1}(x),\\
\omega _{L-1}'(x)&=\omega _{L-1}(x),\\
\tau _{L-1}'(x)&=\tau _{L-1}(x),\\
X_{h,L-1}'(x)&=X_{h,L-1}(x),\  h=1,...,K+1.
\end{align*}

On $[\lf \tau _{L-1}'/\tau _{L-1} \rf\tau _{L-1},\tau _{L-1}')$
we define all the above functions equal to zero with the exception
of the functions $X_{h,L-1}'$,
$h=1,...,K+1$.

When $K>0$
for these functions some minor adjustments
will be made on this interval in order to ensure that
they are pairwise independent  for $h=1,...,K+1$ and are
$M$$-$$0.99$-distributed for $h=1,...,K.$

We can also assume that $X_{K+1,L-1}'$ has constant value
$(1-\rho ')\olll(F_{L-1})\tCC_{\gamma }$
on
$[\lf \tau _{L-1}'/\tau _{L-1} \rf\tau _{L-1},\tau _{L-1}')$.
When $L-1=0$ then this implies that $X_{K+1,L-1}'$
takes this constant value on $\R.$

We define ${E'}^{L-1}$ so that it is periodic by $\tau '_{L-1}$
and
\begin{align}\label{*h16}
&{E'}^{L-1}\cap [0,\tau _{L-1}')=\\
\nonumber
&\left (E^{L-1}\cap \bigg
[0,(\lf \frac{\tau _{L-1}'}{\tau _{L-1}} \rf-1)
\tau _{L-1}\bigg)
\right)\cup \bigg[(\lf \frac{\tau _{L-1}'}{\tau _{L-1}} \rf-1)
\tau _{L-1},\tau _{L-1}'\bigg).
\end{align}
By choosing $\tau _{L-1}'$ sufficiently large we can make
$0\leq {\olll({E'}^{L-1})}-
{\olll(E^{L-1})}$ as small as we wish,
hence, using \eqref{*h10}
there is $\cat_{1}$ such that if
$\tau '_{L-1}>\cat_{1}$ then
\begin{equation}\label{*29x}
\olll({E'}^{L-1})<\frac{L}{(L'+1)}\delta .
\end{equation}

When $L-1=0$ then
\begin{equation}\label{*husv1}
X_{K+1,0}'(x)=(1-\rho ')\olll(F_{0})\tCC_{\gamma }
\text{ for  all }x\in \R.
\end{equation}

If $L-1>0$, that is, $L\geq 2$ we need to deal
with the auxiliary sets related to the distribution
of $X_{K+1,L-1}'.$
We put $S'_{L-1,L-1}=\R$.
Observe that
 \eqref{*h16} holds with
${E'}^{L-1},$ $E^{L-1}$ being replaced by $S'_{L-1,L-1},$
and $S_{L-1,L-1}=\R,$ respectively.
For $l=0,...,L-2$ we define the sets $S'_{L-1,l}$
so that they are periodic by $\tau _{L-1}'$ and we have
$$S'_{L-1,l}\cap [0,\tau _{L-1}')=
S_{L-1,l}\cap [0,(\lf \tau _{L-1}'/\tau _{L-1} \rf-1)\tau _{L-1}).$$
The above definitions and \eqref{*15b} imply
\begin{align}\label{*16b}
X_{K+1,L-1}'(x)&=(1-\rho ')\tCC_{\gamma }\olll(F_{0})
\text{ for }x\in S'_{L-1,0}, \text{ and}\\
\nonumber
X_{K+1,L-1}'(x)&=(1-\rho ')\tCC_{\gamma }\olll(F_{l})
\text{ for }x\in S'_{L-1,l}\sm S'_{L-1,l-1},\  l=1,...,L-1.
\end{align}
By the strict inequalities in \eqref{**h12},
\eqref{*h13} and \eqref{**h13}
we can choose $\cat_{2}$ such that if
$\tau _{L-1}'>\cat_{2}$ then
\begin{equation}\label{*h17}
\frac{1}{1-\rho}
\olll(F_{L-1})>\olll(S_{L-1,0}')>(1-\rho )\olll(F_{L-1})
\end{equation}
and for $l=0,...,L-1$
\begin{equation}\label{**h17}
\frac{1}{(1-\rho )r_{0}\cdot \cdot \cdot r_{l}}
\olll(S_{L-1,0}')>\olll(S_{L-1,l}')>
(1-\rho )\frac{1}{r_{0}\cdot \cdot \cdot r_{l}}
\olll(S_{L-1,0}'),
\end{equation}
or, equivalently,
\begin{equation}\label{***h17}
\frac{1}{1-\rho }\olll(S_{L-1,0}')>
\olll(F_{l})\cdot \olll(S_{L-1,l}')>
(1-\rho )\olll(S_{L-1,0}').
\end{equation}

Set $F'_{L-1}= \{ x:
f'_{K+1,L-1}(x)=1  \}$,
that is, $F'_{L-1}\cap [0,\tau '_{L-1})=
F_{L-1}\cap [0,\lf \tau '_{L-1}/\tau _{L-1} \rf \tau _{L-1})=
F_{L-1}^{\tau '_{L-1}}
\cap [0,\lf \tau '_{L-1}/\tau _{L-1} \rf \tau _{L-1})$ and
$f'_{K+1,L-1}(x)=\chi_{F_{L-1}^{\tau '_{L-1}}}(x)$
for all $x\in\R$.
Clearly, $\olll(F'_{L-1})\leq \olll(F_{L-1})$.

For the case $L-1=0$ we note that $F_{0}\cap [0,\tau _{L-1}')=
[0,\lf \tau '_{L-1}/\tau _{L-1} \rf \tau _{L-1}).$

By \eqref{*14}
for any $x\in \R$ from $\tau '_{L-1}|m$,
it follows that letting $f'_{K+1,L-1}=
\chi_{F_{L-1}^{\tau '_{L-1}}}$ we have
\begin{equation}\label{*f11}
\frac{1}m\sum_{k=n}^{n+m-1}f'_{K+1,L-1}(x+k^{2})
>
(1-\frac{\rho '}2)\olll(F_{L-1}).
\end{equation}
This formula is the main motivation for introducing the
$\tau '_{L-1}$ periodic rearrangements.

Observe that if $\tau _{L-1}'\to\oo$ then
$\olll(F_{L-1}')/\olll(F_{L-1})\to 1$.

Hence we can choose $\cat_{3}$ such that if $\tau '_{L-1}>\cat_{3}$
then
\begin{equation}\label{*85xx}
1-\rho <\frac{\olll(F'_{L-1})}{\olll(F_{L-1})}\leq 1.
\end{equation}
We remark that for $L-1=0$
inequality \eqref{*85xx}
simply means
$1-\rho< \olll(F_{0}') $.

Moreover, we can choose $\cat_{4}$ such that if $\tau '_{L-1}>\cat_{4}$
then
\begin{equation}\label{*19d}
1-\frac{\gamma }{10}<\frac{\olll(F_{L-1}')}{\olll(F_{L-1})}\leq 1
<
1+\frac{\gamma }{10}.
\end{equation}
Finally, by \eqref{**h12} and \eqref{*h17}
if $L\geq 2$
we can choose $\cat_{5}$ such that if $\tau '_{L-1}>\cat_{5}$
then
\begin{equation}\label{*h18}
\frac{1}{1-\rho}\olll(F'_{L-1})>\olll(S'_{L-1,0})>(1-\rho )\olll(F_{L-1}')
\end{equation}
holds as well.

If $x\in [0,\tau _{L-1}')\sm {E'}^{L-1}$, then put
$\alpha _{L-1}'(x)=\alpha_{L-1}(x),$
$\omega _{L-1}'(x)=\omega _{L-1}(x),$
$\tau _{L-1}'(x)=\tau _{L-1}(x)<\tau_{L-1}\ll \tau '_{L-1}.$
This defines $\alpha _{L-1}'(x)$, $\omega _{L-1}'(x)$,
and $\tau _{L-1}'(x)$ for all $x\in \R\sm {E'}^{L-1}$
as well since these functions are periodic by $\tau _{L-1}'$.
It is also clear that $(\tau '_{L-1}(x),p)=1$ for all
$p\in\capp_{L-1}$.
We have
$(\omega _{L-1}'(x))^{2}<\tau _{L-1}$ and
\begin{equation}\label{*19j}
\frac{\omega '_{L-1}(x)}{\alpha _{L-1}'(x)}>\OOO\tau _{L-1}'(x).
\end{equation}
Suppose $\alpha '_{L-1}(x)\leq n<n+m\leq \omega '_{L-1}(x)$
and $\tau '_{L-1}(x)=\tau _{L-1}(x)|m$.
Since  $x\not\in {E'}^{L-1}$,
formula \eqref{*h16} implies that
$x+k^{2}\in [0,\lf \tau _{L-1} '/\tau _{L-1}\rf \tau _{L-1})$
for
$k\in  \{ n,...,n+m-1  \}$
we infer by \eqref{*c17} for $h=1,...,K+1$
\begin{align}\label{*f12}
\frac{1}m\sum_{k=n}^{n+m-1}f_{h,L-1}'(x+k^{2})&=
\frac{1}m\sum_{k=n}^{n+m-1}f_{h,L-1}(x+k^{2})>\\
\nonumber &X_{h,L-1}(x)=X_{h,L-1}'(x).
\end{align}

For all $x\in [0,\tau '_{L-1})\sm {E'}^{L-1}$, $h\in
 \{ 1,...,K+1  \}$,
if ${\alpha'}^{2}_{L-1}(x)\leq j< j+\tau _{L-1}'(x)
\leq {\omega '}^{2}_{L-1}(x),$
$
f_{h,L-1}(x+j+\tau _{L-1}(x))=
f_{h,L-1}'(x+j+\tau '_{L-1}(x))=
f'_{h,L-1}(x+j)=f_{h,L-1}(x+j).$
By periodicity with respect to $\tau '_{L-1}$,
the above estimates hold for any $x\not \in {E'}^{L-1}.$

\subsubsection{Choice of
$\kappa _{L}$,
$q_{L}$, $\FFF_{L}$, $\PPP_{L}$,
and
$X_{h,L}$ on $\R\sm\oLLL'(q_{L})$}\label{425}

Our goal in this section
is to describe some sets, in particular $\FFF_{L}$ and $\PPP_{L}$,
so that we can take $f_{K+1,L} = f_{K+1,L-1}'\chi_{\PPP_{L}}$ in \ref{426}. The functions
$f_{K+1,L-1}'$ and $\chi_{\PPP_{L}}$ are ``independent'' which allows us to
  reduce the integral of $f_{K+1,L}.$ 
In this section we construct three components of  
 the next exceptional set $E^L$. The fourth component will be defined in
 Section \ref{426b}. This last component is the exceptional set coming from
 the $K-M$ family which we put on $\oLLL'(q_{L})$.
We also construct 
 the function $X_{K+1,L}$.   To construct the sets mentioned above
we choose a number $q_{L}\in\capp'_{L},$
$q_{L}=p_{1,L}\cdot \cdot \cdot p_{\kappa _{L},L}$,
$p_{1,L}<...<p_{\kappa_{L} ,L}$,
where $\kappa _{L}$ and $p_{1,L}$ are both sufficiently large.

In fact, we will suppose that $\kappa _{L}$ is larger than the
maximum of the
numbers $\cak_{i}$, $i=1,...,7$ and $\kappa _{\gamma }$,
we also suppose that
$p_{1,L}$ is larger than the maximum of $\pi_i(\kappa_{L})$,
$i=0,...,7$, $p_{1,L}''$, and $p_{\gamma ,\kappa }$,
where $\cak_{i}$, $\pi_{i}$ and
$p_{1,L}''$ are defined below, the numbers $\kappa _{\gamma }$
and $p_{\gamma ,\kappa }$ were defined in Lemma \ref{rdmlem}.
With these assumptions we will be able to use
\eqref{*35x}, \eqref{*44x}, \eqref{*32x}, \eqref{*33x},
\eqref{*21dd}, \eqref{2*33x}, \eqref{*88xxx}, and \eqref{*37yy}
simultaneously.

Recall that we assumed that $\trrr>0$ satisfies \eqref{*38x}.
An application of Lemma \ref{probabilitylemma}
with $\kappa _{L},$ $\epsilon =\delta /4L'$
and
$\trrr$ instead of $\rho$ yields $p_{1,L}''$ sufficiently large so
that $q_{L}=p_{1,L}\cdot \cdot \cdot p_{\kappa_{L},L}$,
with
$p_{1,L}''<p_{1,L}$ satisfies \eqref{*a9} and hence we will be able
to use \eqref{*37yy}.

By \eqref{*9b} for given $\kappa_{L}$ we can choose $\pi_{0}(\kappa _{L})$
such that for $p_{1,L}>\pi_{0}(\kappa _{L})$ we have
\begin{equation}\label{*35x}
\#((\LLL_{0}(q_{L})\sm
\LLL_{0}'(q_{L}))\cap [0,q_{L}))<
\trrr (1-\gamma )\#(\LLL_{0}(q_{L})\cap [0,q_{L})).
\end{equation}

Recall from Remark
\ref{**5dfw}
that the average gap length between points of
$\LLL_{0}(q_{L})$ is approximately
$2^{\kappa _{L}}$ and we can assume that it is
much larger than $\tau _{L-1}'$.
The
normalized
difference between elements of $\LLL_{0}(q_{L})$
approximates
 Poisson distribution by the results in \cite{[Kur]},
see also Lemma \ref{rdmlem}.
We also recall from
Lemma \ref{lputres} that $\oXXX(q_{L})=
\cup_{j\in \Z} [jq_{L} ,jq_{L} +\gamma 2^{\kappa_{L} })$.
We put
\begin{equation}\label{*fi1}
\tFFF_{L}= \{ x:dist(x,\oLLL'(q_{L})\cup
\oXXX(q_{L}))>2\tau '_{L-1}  \},
\end{equation}
\begin{equation}\label{*fi2}
\FFF_{L}=
\cup \{ [j\tau '_{L-1},
(j+1)\tau '_{L-1}):
\tFFF_{L}\cap
[j\tau '_{L-1},(j+1)\tau '_{L-1})\not=\ess
\},
\end{equation}
$$\hFFF_{L}= \{ x:dist(x,\oLLL'(q_{L}))>\tau '_{L-1}  \},$$
$$\hPPP_{L}= \{ x:dist(x,\R\sm\oLLL'(q_{L}))\leq\tau '_{L-1}  \},$$
$$
\PPP_{L}=
\cup \{ [j\tau '_{L-1},
(j+1)\tau '_{L-1}):
\hPPP_{L}\cap
[j\tau '_{L-1},(j+1)\tau '_{L-1})\not=\ess
\},
$$
and finally
$$\tPPP_{L}= \{ x:dist(x,\R\sm\oLLL'(q_{L}))\leq 2\tau '_{L-1}  \}.$$
It is clear that
\begin{equation}\label{*17b}
\tFFF_{L}\sse\FFF_{L}\sse \hFFF_{L}\sse \R\sm \oLLL'(q_{L})
\sse \hPPP_{L}\sse \PPP_{L}
\sse \tPPP_{L}.
\end{equation}
It is also important that by \eqref{*fi1} and \eqref{*fi2} we have
\begin{equation}\label{*fi3}
\FFF_{L}\cap \oXXX(q_{L})=\ess.
\end{equation}
The sets
$\FFF_{L}$ and
$\PPP_{L}$ are periodic by $\tau '_{L-1} q_{L}$
and the sets $\tFFF_{L}$, $\hFFF_{L}$, $\tPPP_{L}$,
and $\hPPP_{L}$ are periodic by $q_{L}$.
If $\kappa _{L}$ is sufficiently large, then
$2^{\kappa _{L}}$ and hence most of the gaps between points
of $\LLL_{0}(q_{L})$ are much larger than $\tau '_{L-1}.$
In the sequel by $\approx$ we mean that if $p_{1,L}$
and $2^{\kappa _{L}}$ (compared to $\tau'_{L-1}$) are
sufficiently large, then the ratio of the two sides of
$\approx$ is sufficiently close to $1$, later
we will specify further this assumption.
Since $\oLLL'(q_{L})$ consists of intervals of length
$\gamma 2^{\kappa_{L} }$
which is much larger than $\tau _{L-1}'$,
 we have
\begin{align}
\label{*16}
&\olll(\tPPP_{L})\approx\olll(\tFFF_{L})
\approx\olll(\PPP_{L})\approx\olll(\FFF_{L})\approx
\olll(\hPPP_{L})\approx\olll(\hFFF_{L}),\\
\label{**16}
&\olll(\FFF_{L})\leq\olll( \R\sm \oLLL'(q_{L}))\leq
\olll( \PPP_{L})
\text{ and }
\olll(\R\sm \oLLL'(q_{L}))\approx \olll(\PPP_{L})
.
\end{align}

Using this and \eqref{*19bud} we can choose $\cak_{1}$ and a function
$\pi_{1}$ such that if $\kappa _{L}>\cak_{1}$
and $p_{1,L}>\pi_{1}(\kappa _{L})$ then
\begin{equation}\label{*44x}
\olll(\FFF_{L})>
\olll(\R\sm\oLLL'(q_{L}))/2>(1-\gamma )/2.
\end{equation}

Set $E_{L}''=\R\sm (\oLLL'(q_{L})\cup \tFFF_{L}),$
this will be part of the new exceptional set
$E^{L}$.
We also introduce 
\begin{equation}\label{*080211c}
\tEE_{L}''=
\R\sm (\oLLL'(q_{L})\cup \FFF_{L})
\sse E_{L}''.\end{equation}
 It is clear that
 $E_{L}''$ is periodic by $q_{L}$, while $\tEE_{L}''$ is periodic
by $\tau _{L-1}' q_{L}$.

We can choose $\cak_{2}$ and a function
$\pi_{2}$ such that if $\kappa _{L}>\cak_{2}$
and $p_{1,L}>\pi_{2}(\kappa _{L})$ then
\begin{equation}\label{*32x}
\olll(\tEE_{L}'')
\leq\olll(E''_{L})<
\frac{\delta}{4L'}.
\end{equation}

Set
\begin{equation}\label{*18b}
X_{h,L}(x)=X_{h,L-1}'(x)\text{ if }h\leq K+1\text{ and }
x\in \FFF_{L}.
\end{equation}

For $K>0$ we can make the following comment:
Since $\FFF_{L}$ consists of intervals of the form
$[j\tau '_{L-1},(j+1)\tau '_{L-1}),$ this definition
and the remark after the definition of
$X_{h,L-1}'$ in Section \ref{424} ensures that the functions
$X_{h,L}$ are pairwise independent  for $h=1,...,K+1$ and are
conditionally $M$$-$$0.99$ distributed on $\FFF_{L}$
for $h=1,...,K$.

On $\tEE''_{L}$
we will have
\begin{equation}\label{*080211b}
X_{K+1,L}(x)=(1-\rho ')\olll(F_{L-1})\tCC_{\gamma }
\end{equation}
and
 we define $X_{h,L}$ for $h=1,...,K$
so that they are pairwise independent  on $\tEE''_{L}$,
furthermore, (for $K>0$) the functions
$X_{h,L}$ are conditionally $M$$-$$0.99$-distributed on $\tEE''_{L}$
for $h=1,...,K.$ Since $X_{K+1,L}$ is constant on $\tEE''_{L}$
it is automatically independent on this set from $X_{h,L}$ for $h=1,...,K.$ 

The functions
$X_{h,L}$ are periodic on $\tEE_{L}''$ by $\tau '_{L-1}q_{L}$
for $h=1,...,K+1.$ In this way the $X_{h,L}$'s are defined
on $\FFF_{L}\cup \tEE_{L}''=\R\sm \oLLL'(q_{L})$.

We set $E_{L}'=\PPP_{L}\cap {E'}^{L-1}.$

Next we consider some sets which are used to
describe the distribution of $X_{K+1,L}$.

If $L=1$ set $S_{1,0}=\R\sm \oLLL'(q_{1})$ and
$S_{1,1}=\R.$

If $L\geq 2$
first we define $S_{L,l}$
for $l=0,...,L-1$ so that
$S_{L,l}\cap \FFF_{L}=S'_{L-1,l}\cap \FFF_{L}$
for $l=0,...,L-1.$
We choose $S_{L,l}$ so that $S_{L,l}\cap (\R\sm \FFF_{L})=
\ess$
for $l=0,...,L-2.$
We choose $S_{L,L-1}$ so that
$S_{L,L-1}=\tEE''_{L}\cup (S'_{L-1,L-1}\cap \FFF_{L})=
\tEE''_{L}\cup\FFF_L=\R\sm \oLLL'(q_L).$
Finally, we set $S_{L,L}=\R,$
then
$S_{L,L}\sm S_{L,L-1}=
\oLLL'(q_{L}).$

We have by \eqref{*16b},
(for the case
$L=1$ by \eqref{*husv1})
 and \eqref{*18b}
\begin{align}\label{**18b}
& X_{K+1,L}(x)=(1-\rho ')\olll(F_{0})\tCC_{\gamma }
\text{ for }x\in S_{L,0}, \text{ and }\\
\nonumber
& X_{K+1,L}(x)=(1-\rho ')\olll(F_{l})\tCC_{\gamma }
\text{ for }x\in S_{L,l}\sm S_{L,l-1},
\  l=1,...,L-1.\end{align}
The case when $l=L$ will be considered in
\eqref{*21d}.
By \eqref{*080211c} we have $\tEE''_{L}\cap \FFF_{L}=\ess$
and hence $S_{L,l}\cap \tEE''_{L}=\ess$ for $l\leq L-2$.
This implies $\tEE''_{L}\sse S_{L,L-1}\sm S_{L,L-2}$.
Hence  \eqref{**18b} applied with $l=L-1$ implies
\eqref{*080211b}.

Let $F_{L}=\PPP_{L}\cap F'_{L-1}.$
Using the fact that $F'_{L-1}$ is periodic
by $\tau _{L-1}'$ and $\PPP_{L}$ is the union
of some intervals of the form
$[j\tau '_{L-1}, (j+1)\tau '_{L-1})$
and is periodic by $\tau '_{L-1}q_{L}$
one can easily see that
$\olll(F_{L})=\olll(\PPP_{L})\olll(F_{L-1}').$
Moreover,
$$r_{L}\defeq \frac{\olll(F_{L})}
{\olll(F_{L-1})}
=\olll(\PPP_{L})\frac{\olll(F'_{L-1})}
{\olll(F_{L-1})}\approx \olll(\PPP_{L})\approx\olll(\FFF_{L})
\approx \olll(\R\sm\oLLL'(q_{L})).$$
Since $\olll(F'_{L-1})\leq \olll(F_{L-1})$
by \eqref{**16} there is
$\cak_{3}$ and a function
$\pi_{3}$ such that if $\kappa _{L}>\cak_{3}$
and $p_{1,L}>\pi_{3}(\kappa _{L})$ then
\begin{equation}\label{*33x}
\olll(\PPP_{L})\frac{\olll(F_{L-1}')}
{\olll(F_{L-1})}<\frac{1-\gamma }
{1-\gamma -\gamma ^2}\olll(\R\sm\oLLL'(q_{L})).
\end{equation}
By \eqref{*13/1} we obtain
\begin{equation}\label{*21j}
\tCC_{\gamma }r_{L}=\tCC_{\gamma }\frac{\olll(F_{L})}
{\olll(F_{L-1})}=\tCC_{\gamma }\olll(\PPP_{L})
\frac{\olll(F_{L-1}')}
{\olll(F_{L-1})}<1-\gamma .
\end{equation}
By \eqref{*19bud},
\eqref{*15d}, \eqref{*19d} and $\olll(\PPP_{L})\approx
\olll(\R\sm \oLLL'(q_{L}))$
there is $\cak_{4}$ and a function
$\pi_{4}$ such that if $\kappa _{L}>\cak_{4}$
and $p_{1,L}>\pi_{4}(\kappa _{L})$ then
\begin{equation}\label{*21dd}
1-2\gamma <r_{L}= \frac{\olll(F_{L})}{\olll(F_{L-1})}
=\olll(\PPP_{L})\frac{\olll(F'_{L-1})}{\olll(F_{L-1})}
<1-\frac{\gamma }2.
\end{equation}
We set
\begin{equation}\label{*21d}
X_{K+1,L}(x)=(1-\rho ')\olll(F_{L})\tCC_{\gamma }
\text{ for }x\in S_{L,L}\sm S_{L,L-1}=\oLLL'(q_{L}).
\end{equation}

If $L=1$ then $\olll(F_{1})=\olll(\PPP_{1})\olll(F_{0}')\approx
\olll(\R\sm \oLLL'(q_{1}))$ and $\olll(S_{1,0})=\olll(\R\sm
\oLLL'(q_{1})).$ Furthermore, $\olll(S_{1,1})=\olll(\R)=1,$
$r_{1}\approx \olll(\R\sm \oLLL'(q_{1}))$.

By \eqref{*19bud}, \eqref{*13/1}, \eqref{*16} and \eqref{**16}
there is
$\cak_{5}$ and a function
$\pi_{5}$ such that if $\kappa _{1}>\cak_{5}$
and $p_{1,1}>\pi_{5}(\kappa _{1})$ then
\begin{equation}\label{2*33x}
1-1.1\gamma <\olll(\PPP_{1})<1-\frac{8\gamma }{10}\text{ and }
1-\rho < \frac{\olll(\R\sm\oLLL'(q_{1}))}{\olll(\PPP_{1})}\leq 1.
\end{equation}

From \eqref{*85xx}, \eqref{*19d}, \eqref{**16} and \eqref{2*33x} it follows that
\begin{equation}\label{50401}
1-2\gamma< r_1=\olll(\PPP_{1})\frac{\olll(F_{0}')}{\olll(F_{0})}
< 1-\frac{\gamma}2,\end{equation}
$$
\frac{1}{1-\rho }\olll(F_{1})=
\frac{1}{1-\rho }\olll(\PPP_{1})\olll(F_{0}')
>$$
$$\olll(S_{1,0})
=
\olll(\R\sm \oLLL'(q_1))>
(1-\rho )\olll(\PPP_{1})\olll(F_{0}')= (1-\rho )\olll(F_{1})
,$$
(keeping in mind $F_{0}=\R$)
\begin{equation}\label{87a}
\frac{1}{(1-\rho )r_{0}r_{1}}
\olll(S_{1,0})=\frac{\olll(\R\sm \oLLL'(q_{1}))\olll(F_{0})}
{(1-\rho )\olll(\PPP_{1})\olll(F_{0}')}>\olll(S_{1,1})=1>
\end{equation}
$$
(1-\rho )
\frac{\olll(\R\sm \oLLL'(q_{1}))\olll(F_{0})}
{\olll(\PPP_{1})\olll(F_{0}')}=
(1-\rho )\frac1{r_{0}
r_{1}}\olll(S_{1,0})
$$
and
\begin{equation}\label{88a}
\frac1{1-\rho }\olll(S_{1,0})>\olll(F_{1})
\olll(S_{1,1})>(1-\rho )\olll(S_{1,0}).
\end{equation}
This shows that \eqref{*h26} and \eqref{**h26}
below
hold for $L=1$ and $l=1.$
For $L=1$ and $l=0$, \eqref{*h26} and \eqref{**h26}
are obvious.

If $L\geq 2$
we have $\olll(S_{L,0})=\olll(\FFF_{L})\cdot \olll(S_{L-1,0}')$.
From \eqref{*h18} it follows that
\begin{equation}\label{*1rho}
1>(1-\rho )\frac{\olll(F'_{L-1})}{\olll(S'_{L-1,0})}.
\end{equation}
Therefore, by \eqref{*16} and \eqref{**16}
there is
$\cak_{6}$ and a function
$\pi_{6}$ such that if $\kappa _{L}>\cak_{6}$
and $p_{1,L}>\pi_{6}(\kappa _{L})$ then
\begin{equation}\label{*88xxx}
1\geq \frac{\olll(\FFF_{L})}
{\olll(\PPP_{L})}
>(1-\rho )\frac{\olll(F'_{L-1})}{\olll(S'_{L-1,0})}.
\end{equation}
Using this, \eqref{*h18} and \eqref{**16} a simple calculation shows that
\begin{equation}\label{*h25}
\frac{1}{1-\rho}\olll(F_{L})
=\frac{1}{1-\rho}\olll(\PPP_{L})\olll(F_{L-1}')
>\olll(S_{L,0})=\olll(\FFF_{L})\olll(S_{L-1,0}')>
\end{equation}
$$(1-\rho )\olll(\PPP_{L})\olll(F_{L-1}')=
(1-\rho )\olll(F_{L}).$$
It is also clear that $\olll(S_{L,l})=\olll(\FFF_{L})\cdot
\olll(S_{L-1,l}')$ for $l=0,...,L-2$.
Using \eqref{**h17}
 we have
for $l=0,...,L-2$
\begin{equation}\label{*h26}
\frac{1}{(1-\rho )r_{0}\cdot \cdot \cdot r_{l}}
\olll(S_{L,0})>\olll(S_{L,l})>(1-\rho )\frac1{r_{0}\cdot \cdot \cdot
r_{l}}\olll(S_{L,0})
\end{equation}
and by $\olll(F_{l})=r_{0}\cdot \cdot \cdot r_{l}$ we have
\begin{equation}\label{**h26}
\frac1{1-\rho }\olll(S_{L,0})>\olll(F_{l})
\olll(S_{L,l})>(1-\rho )\olll(S_{L,0}).
\end{equation}

When $l=L-1$ a little caution is needed.
We have $\olll(S_{L,L-1})=(\olll(\FFF_{L})+\olll(\tEE_{L}''))
\olll(S'_{L-1,L-1})=\olll(\R\sm \oLLL'(q_{L}))$.
By \eqref{**h17}
$$\frac{1}{(1-\rho )r_{0}\cdot \cdot \cdot r_{L-1}}\cdot
\frac{\olll(S'_{L-1,0})}{\olll(S'_{L-1,L-1})}>1.$$
Hence,
there is
$\cak_{7}\geq \cak_{1}$ and a function
$\pi_{7}\geq \pi_{1}$ such that if $\kappa _{L}>\cak_{7}$
and $p_{1,L}>\pi_{7}(\kappa _{L})$ then using \eqref{*44x}
and \eqref{*32x}
\begin{equation}\label{*44xx}
\frac{1}{(1-\rho )r_{0}\cdot \cdot \cdot r_{L-1}}\cdot
\frac{\olll(S'_{L-1,0})}{\olll(S'_{L-1,L-1})}>
{1+\frac{2\olll(E_{L}'')}{1-\gamma }}
>\end{equation}
$$
{1+\frac{\olll(\tEE''_{L})}{\olll(\FFF_{L})}}=
\frac{\oLLL(\FFF_{L})+\olll(\tEE_{L}'')}{\olll(\FFF_{L})}.
$$
Using this and \eqref{**h17} one can deduce
that \eqref{*h26} and \eqref{**h26} hold when
$l=L-1$.

From $\olll(S_{L,L})=1$ and \eqref{*h25} it follows that
$$\frac1{1-\rho }\olll(S_{L,0})>
\olll(F_{L})\olll(S_{L,L})>(1-\rho)\olll(S_{L,0}).$$
Using the fact that $\olll(F_{L})=r_{0}\cdot \cdot \cdot r_{L}$
we find that \eqref{*h26} and \eqref{**h26}
hold for $l=L$ as well.

Denote by $\oEE'''_{L}$ the set of those $n$'s for which
\begin{equation}\label{*EvvL}
\# \bigg(((n+\LLL_{0}(q_{L}))\sm \LLL(q_{L}))\cap [0,q_{L})\bigg)<
(1-\trrr)(1-\gamma )\#(\LLL_{0}(q_{L})\cap [0,q_{L})).
\end{equation}
Recall that by our choice $p_{1,L}>p'_{1,L}$ and hence
$q_{L}$ satisfies \eqref{*a9} with $\epsilon =\delta /4L'$
and $\trrr$ instead of $\rho$.
Lemma \ref{probabilitylemma} yields
\begin{equation}\label{*37yy}
\#(\oEE_{L}'''\cap [0,q_{L}))<\frac{\delta }{4L'}
q_{L}.
\end{equation}
Set 
\begin{equation}\label{*EvvL2}
E_{L}'''= \{ x\in \R:
\lf x \rf\in \oEE'''_{L}  \}.
\end{equation}
Then $\olll(E_{L}''')<\delta /4L'$ and if
$x\not\in E_{L}'''$ we have
\begin{equation}\label{*h28}
\frac{1}{q_{L}}\# \{ k'\in
[0,q_{L})\cap \Z
:\lf x \rf+{k'}^{2}\not\in\LLL(q_{L})  \}\geq
\end{equation}
$$\frac{1}{q_{L}}\# \{ k'\in
[0,q_{L})\cap \Z
:\lf x \rf+{k'}^{2}\not\in\LLL(q_{L}),
{k'}^{2}\in \LLL_{0}'(q_{L})  \}\geq$$
$$\frac1{q_{L}}2^{\kappa _{L}}\#
(((\lf x \rf+\LLL_{0}'(q_{L}))\sm \LLL(q_{L}))\cap[0,q_{L}))\geq$$
$$\frac1{q_{L}}2^{\kappa _{L}}
\bigg(
\#
(((\lf x \rf+\LLL_{0}(q_{L}))\sm \LLL(q_{L}))\cap[0,q_{L}))-
\#
((\LLL_{0}(q_{L})\sm\LLL_{0}'(q_{L}))\cap
[0,q_{L}))\bigg)\geq$$
(using that for $\lf x \rf\not \in \oEE_{L}'''$ we have the negation of \eqref{*EvvL})
$$\frac{1}{q_{L}}2^{\kappa _{L}}\bigg(
((1-\trrr)(1-\gamma )
\#(\LLL_{0}(q_{L})\cap
[0,q_{L})))-\#
((\LLL_{0}(q_{L})\sm \LLL_{0}'(q_{L}))\cap
[0,q_{L}))\bigg)=(*).$$
Recall that we assumed that $p_{1,L}>\pi_{0}(\kappa _{L})$.
Thus we can apply \eqref{*35x}
yielding
$$(*)\geq \frac{1}{q_{L}}2^{\kappa _{L}}
(1-2\trrr)(1-\gamma )\#(\LLL_{0}(q_{L})\cap [0,q_{L}))=(**).$$
Using \eqref{*3bud}
we can finish with the inequality
\begin{equation}\label{*h29}(**)>
\frac{1}{q_{L}}2^{\kappa _{L}}(1-2\trrr)(1-\gamma )
\frac{q_{L}}{2^{\kappa _{L}}}=
(1-2\trrr)(1-\gamma ).\end{equation}

\subsubsection{Putting $K-M$ families on $\oLLL'(q_{L})$}
\label{426}

In this section we check the domination property of averages
along squares for one part of the complement of $E^L.$

We put
\begin{equation}\label{*20b}
f_{K+1,L}=f'_{K+1,L-1}\cdot\chi _{\PPP_{L}}.
\end{equation}
Then indeed, $F_{L}=
 \{ x: f_{K+1,L}(x)=1  \}=\PPP_{L}\cap F_{L-1}',$
and we have $\ointt f_{K+1,L}=\olll(F_{L}).$
Since $f'_{K+1,L-1}$ is periodic by $\tau '_{L-1}\in \capp'_{L-1}$
and $\PPP_{L}$ is periodic by $\tau '_{L-1}q_{L}$
the function $f_{K+1,L}$ is
also periodic by $\tau '_{L-1}q_{L}$.

When $K>0$ and $h\in \{ 1,...,K \}$ we will define $f_{h,L}$ in \eqref{*19}
so that $$f_{h,L}(x)=f'_{h,L-1}(x)\text{ for }x\in \FFF_{L}.$$

Choose $\capp''_{L}\sse \capp'_{L}$ such that
it contains infinitely many primes and all
their possible products, moreover
all numbers in
$\capp_{L}''$ are relatively prime to $q_{L}\in \capp'_{L}$
and set $\ocap_{L}=\capp\cup \capp_{L}''\cup
\capp_{L+1}'\cup...\cup \capp'_{L'}.$

When $K=0$, {\it i.e.}, when constructing a $1-M$ family we have
to put a vacuous ``$0-M$ family" on $\oLLL'(q_{L})$. 

Hence, for $K=0$, for the definition of a 
$1-M$ family we just set $\ottt_{L}=q_{L}
(\tau '_{L-1})^{3}$,
$E_{\delta _{L}}=\ess$. 

When $K\geq 1$
this step is crucial, see \ref{426b}. 

Set $\capp_{L}=\capp\cup \capp_{L+1}'\cup
...\cup\capp_{L'}'\sse \ocap_{L} \sse \capp_{L-1}$.
We also put
$E^{L}=E_{\delta _L}\cup E_{L}'\cup E_{L}''\cup E_{L}'''$,
and
 $\tau _{L}=q_{L}\ottt_{L}\tau _{L-1}'$.
When $K=0$ we defined $\ottt_{L}=q_{L}
(\tau '_{L-1})^{3}$ and hence  $\tttt_{L}=q_{L}^2(\tau _{L-1}')^4$.
Then
for all $p\in \capp_{L}\sse \ocap_{L}$
we have $(\tau _{L},p)=1$ and $(\tau'_{L-1}(x),p )=1$ when $x\in \R\sm E^{L}$.
Assume $x\in \R\sm(\oLLL'(q_{L})\cup E^{L}).$
Then
$x\in \tFFF_{L}\sse \FFF_{L}$ and
 the old estimates work.
 
In other words, for $x\in \R\sm (\oLLL'(q_{L})\cup E^{L})
\sse \tFFF_{L},$
set $\alpha _{L}(x)=\alpha _{L-1}'(x),$
$\omega _{L}(x)=\omega '_{L-1}(x),$
$\tau _{L}(x)=\tau '_{L-1}(x)=\tau _{L-1}(x)$.
Then
$\omega _{L}^{2}(x)<\tau _{L-1}<\tttt'_{L-1}<\tau _{L}$
and by \eqref{*19j} we have
$$
\frac{\omega _{L-1}'(x)}
{\alpha _{L-1}'(x)}=
\frac{\omega _{L}(x)}{\alpha _{L}(x)}>
\OOO \tau _{L-1}'(x)=\OOO \tau _{L}(x).$$

Observe that if $x\in \R\sm (\oLLL'(q_{L})\cup E^{L})
\sse \tFFF_{L}\sm E'_{L}=\tFFF_{L}\sm {E'}^{L-1}
\sse \tFFF_{L}\sse \FFF_{L}$,
then from $x\not\in {E'}^{L-1},$
 $\alpha _{L}(x)\leq n\leq k<n+m\leq \omega _{L}(x)$,
and
$\omega ^{2}_{L}(x)<\tau _{L-1}<\tau'_{L-1}$ it follows that
$x+k^{2}\in \FFF_L\sse \PPP_{L}$ and hence by \eqref{*20b},
$f_{K+1,L}(x+k^{2})=f'_{K+1,L-1}(x+k^2)$
and by \eqref{*f12} , if $\tau _{L}(x)=\tau'_{L-1}(x)|m$ then
$$\frac1{m}\sum_{k=n}^{n+m-1}f'_{K+1,L-1}(x+k^{2})
=\frac1{m}\sum_{k=n}^{n+m-1}f_{K+1,L}(x+k^{2})
>
$$ $$
X_{K+1,L}(x)=X_{K+1,L-1}'(x).$$

Finally, for all $x\in \R\sm(\oLLL'(q_{L})\cup E^{L})
\sse \tFFF_{L}\sm E'_{L}=\tFFF_{L}\sm {E'}^{L-1}
\sse \tFFF_{L}$,
$h=1,...,K+1$, if
$\alpha _{L}^{2}(x)\leq j<j+\tau _{L}(x)\leq \omega _{L}^{2}(x)
<\tau_{L-1}$,
then $x+j,\, x+j+\tttt_{L}(x)\in \FFF_{L}$ and
$$f'_{h,L-1}(x+j+\tau '_{L-1}(x))=
f_{h,L}(x+j+\tau _{L}(x))=f_{h,L}(x+j)=
f'_{h,L-1}(x+j).$$

$\dagger$ Assume $x\in \oLLL'(q_{L})\sm E^{L}\sse\oLLL'(q_{L})\sm
E_{\delta _L}$.

For the $K=0$ case set $\alpha _{L}(x)=\alpha _{L-1}'(x),$
$\omega _{L}(x)=\tau '_{L-1}q_{L}\omega '_{L-1}(x),$  
 and
$\tau _{L}(x)=\tau '_{L-1}q_{L}\tau '_{L-1}(x)$.
Then $\omega _{L}^{2}(x)<(\tau '_{L-1})^{2}q_{L}^{2}(\tau '_{L-1})^{2}=
\tau _{L}$ and $$\frac{\omega _{L}(x)}{\alpha _{L}(x)
}=\frac{\omega '_{L-1}(x)\tau '_{L-1}q_{L}}{\alpha '_{L-1}(x)}>
\OOO \tau '_{L-1}(x)\tau '_{L-1}q_{L}=\OOO \tau _{L}(x).$$
Since $f_{K+1,L}$ is periodic by $\tau '_{L-1}q_{L}$
and $\tau '_{L-1}q_{L}|\tau _{L}(x)$ we have
$f_{K+1,L}(x+j+\tau _{L}(x))=f_{K+1,L}(x+j)$
for all $x$ and $j$.

 Instead of the above
paragraph we will have a different argument in Subsection \ref{426b}
for the $K>0$ case.

\subsubsection{Properties of $f_{K+1,L}$}
\label{427}

In this section we check the domination property for averages
along squares for $x$ in the remaining part of the complement of $E^L.$

We need to check \eqref{*c17} when $h=K+1$ and
$x\in \oLLL'(q_{L})\sm E^{L}=
(S_{L,L}\sm S_{L,L-1})\sm E^{L},$
and
$\tau _{L}(x)
=\tau '_{L-1}q_{L}\tau'_{L-1}(x)|m$.
If we can show that \eqref{*c17}
holds
 when $m=\tau '_{L-1}q_{L}$
then this clearly implies that it holds when
$\tau _{L}(x)| m$.
Recall that $\hPPP_{L}$ is periodic by $q_{L}.$
Since
$\tau'_{L-1}\in \capp_{L-1}'$, $q_{L}\in \capp_{L}'$
implies
$(\tau '_{L-1},q_{L})=1$, $k'+jq_{L}$ covers all
residues modulo $\tau _{L-1}'$ as $j$ runs from $0$
to $\tau _{L-1}'-1.$
Since $f'_{K+1,L-1}$
is periodic by $\tau _{L-1}'$, using \eqref{*f11}
we obtain
\begin{equation}\label{*22b}
\frac1{\tau '_{L-1}}\sum _{j=0}^{\tau _{L-1}'-1}f'_{K+1,L-1}
(x+(k'+jq_{L})^{2})> (1-\frac{\rho '}2)\olll(F_{L-1}).
\end{equation}
Also observe that from the periodicity of $\hPPP_{L}$
by $q_{L}$ it follows that if $x+{k'}^{2}\in\hPPP_{L}$ then
$x+(k'+jq_{L})^{2}\in \hPPP_{L}\sse \PPP_{L}$ as well.
Hence from $x+{k'}^2\in\hPPP_{L}$ it follows that
$f_{K+1,L}(x+(k'+jq_{L})^{2})=f'_{K+1,L-1}
(x+(k'+jq_{L})^{2}).$
Therefore,
$$\frac1{\tau '_{L-1}q_{L}}
\sum _{k=n}^{n+\tau _{L-1}'q_{L}-1}f_{K+1,L}
(x+k^{2})=$$ $$
\frac1{\tau '_{L-1}q_{L}}
\sum_{k'=n}^{n+q_{L}-1}\sum _{j=0}^{\tau _{L-1}'-1}f_{K+1,L}
(x+(k'+jq_{L})^{2})\geq
$$
$$\frac1{\tau '_{L-1}q_{L}}
\sum_{
\begin{array}{c}
\scriptstyle{k'=n}\\
\scriptstyle{x+{k'}^2\in\hPPP_{L}}
\end{array}}^{n+q_{L}-1}\sum _{j=0}^{\tau _{L-1}'-1}f_{K+1,L}
(x+(k'+jq_{L})^{2})=$$
$$\frac1{q_{L}}
\sum_{
\begin{array}{c}
\scriptstyle{k'=n}\\
\scriptstyle{x+{k'}^2\in\hPPP_{L}}
\end{array}}^{n+q_{L}-1}
\frac{1}{\tau '_{L-1}}
\sum _{j=0}^{\tau _{L-1}'-1}f'_{K+1,L-1}
(x+(k'+jq_{L})^{2})>$$
(using \eqref{*22b})
$$
\frac1{q_{L}}
\sum_{
\begin{array}{c}
\scriptstyle{k'=n}\\
\scriptstyle{x+{k'}^2\in\hPPP_{L}}
\end{array}}^{n+q_{L}-1}
(1-\frac{\rho '}2)\olll(F_{L-1})\geq
$$
(using \eqref{*17b})
$$
\frac{1}{q_{L}}
(1-\frac{\rho '}2)\olll(F_{L-1})
\#  \{ k'\in [0,q_{L})\cap \Z: x+{k'}^{2}\not\in
\oLLL'(q_{L})  \}\geq
$$
$$
\frac{1}{q_{L}}
(1-\frac{\rho '}2)\olll(F_{L-1})
\#  \{ k'\in [0,q_{L})\cap \Z: x+{k'}^{2}\not\in
\oLLL(q_{L})  \}=
$$
$$
(1-\frac{\rho '}2)\olll(F_{L-1})
\frac{1}{q_{L}}\#  \{ k'\in [0,q_{L})\cap \Z: \lf x \rf+{k'}^{2}\not\in
\LLL(q_{L})  \}.
$$
Now use the estimates \eqref{*h28} through \eqref{*h29}
and obtain that for $x\not\in E_{L}'''\sse E^{L}$
$$\frac1{q_{L}}
\# \{ k'\in [0,q_{L})\cap \Z:
\lf x \rf+{k'}^{2}\not\in
\LLL(q_{L})  \}\geq
(1-2\trrr)(1-\gamma ).$$
Thus, if $x\in (S_{L,L}\sm S_{L,L-1})\sm E^{L}=\oLLL'(q_{L})\sm E^{L}$ we have
$$\frac1{\tau '_{L-1}q_{L}}
\sum _{k=n}^{n+\tau _{L-1}'q_{L}-1}f_{K+1,L}
(x+k^{2})>
$$ $$
(1-\frac{\rho '}2)
\olll(F_{L-1})(1-2\trrr)
(1-\gamma )\geq$$
(Using \eqref{*38x}, \eqref{*21j} and \eqref{*21d})
$$\geq (1-\rho ')(1-\gamma )\olll(F_{L-1})>
(1-\rho ')\tCC_{\gamma }\olll(F_{L})=X_{K+1,L}(x).$$

\subsubsection{Finishing the leakage}\label{428}

We keep repeating the leakage steps until for the first time
for  some $L''$ we have $\olll(F_{L''})<2^{-M}$ which implies 
 $\olll(F_{L''-1})\geq 2^{-M}$.
By \eqref{*h6} and \eqref{*21dd} applied  to all $L\leq L'$
we have $L''\leq L'$
and by $\gamma <\gamma _{0}<10^{-7}$ we have $L''\geq 2.$

We set $f_{h}=f_{h,L''}$ for $h=1,...,K+1,$
and $X_{h}=X_{h,L''}$ for $h=1,...,K.$
From the induction steps we have
$E_{\delta }\defeq E^{L''}$ such that $\olll(E_{\delta })
< \frac{(L''+1)}{(L'+1)}\delta \leq \delta .$
There exists $\tau \defeq \tau _{L''}$ such that
$f_{h}$, $h=1,...,K+1$, $X_{h}$, $h=1,...,K$
and $X_{K+1,L''}$ are periodic by $\tau ,$
$X_{h}$, $h=1,...,K$ and $X_{K+1,L''}$ are pairwise independent
$X_{h}$, $h=1,...,K$ are $M$$-$$0.99$-distributed.
By using
the distributional properties of
$X_{K+1,L''}$ we will define
$X_{K+1}$ at the end of this section.

For all $x\not\in E_{\delta }$ there exist
$\omega (x)=\omega _{L''}(x)>\alpha (x)
=\alpha _{L''}(x)>A,$ $\tau (x)=\tau _{L''}(x)
< \tau $
such that
$\omega ^{2}(x)<\tau ,$
$\frac{\omega (x)}{\alpha (x)}>\OOO \tau (x)$.
Setting $f_{K+1}=f_{K+1,L''}$,
(see also \eqref{*20b}) if $\tau (x)|m$ then
\begin{equation}\label{*24b}
\frac{1}m\sum_{k=n}^{n+m-1}f_{K+1}(x+k^{2})>X_{K+1,L''}(x).
\end{equation}
When $K>0$ one also needs to use \eqref{*h44}, see Section \ref{428b}.

For all $p\in\capp_{L''}\supset\capp$,
$(\tau (x),p)=1,$ $(\tau ,p)=1.$
For all $x\not\in E_{\delta }$
and for all $h\in  \{ 1,...,K +1 \}$,
$f_{h}(x+j+\tau (x))=f_{h}(x+j)$ whenever
$\alpha ^{2}(x)\leq j<j+\tau (x)\leq \omega ^{2}(x).$
Finally,
\begin{equation}\label{*14bx}
\ointt f_{K+1}=\olll(F_{L''})<2^{-M+1}<
\GGG\cdot 2^{-M+1},
\end{equation}
when $K>0$ we also need \eqref{*48yy} from Section \ref{428b}.

We have met all the requirements for a $(K+1)-M$ family except the
distribution of $X_{K+1,L''}$ is not quite right.
We need to replace $X_{K+1,L''}$ by a suitably chosen
$X_{K+1}$ which is $M$$-$$0.99$-distributed,
moreover for $K>0$ it is
pairwise independent from $X_{h}$ when $h=1,...,K$.
By choosing
$X_{K+1}$ so that $X_{K+1}\leq X_{K+1,L''}$
from \eqref{*24b} we infer that
if $\alpha (x)\leq n<n+m\leq \omega (x)$
and
$\tau (x)|m$ then
\begin{equation}\label{*h44bb}
\frac{1}m\sum_{k=n}^{n+m-1}f_{K+1}(x+k^{2})>X_{K+1}(x).
\end{equation}

Since $L''$ is the first index when
$\olll(F_{L''})<2^{-M}$
we have $\olll(F_{L''-1})\geq 2^{-M}$ which
by \eqref{*21dd}
implies
\begin{equation}\label{*28j}
\olll(F_{L''})>(1-2\gamma )2^{-M}.
\end{equation}

\medskip
What is the distribution of $X_{K+1,L''}$?
Recall $F_{L''}= \{ x:f_{K+1}(x)=1  \}$,
$1-2\gamma <r_{L}=\olll(F_{L})/\olll(F_{L-1})<1-\frac{\gamma}2 ,$
for $L=1,...,L'',$ and
$\olll(F_{L})=r_{0}\cdot \cdot \cdot r_{L}=
r_{1}\cdot \cdot \cdot r_{L}$.
By \eqref{*h25} and \eqref{*39x}
\begin{equation}\label{*27c}
\frac{1}{1-\rho}2^{-M}>
\frac{1}{1-\rho}\olll(F_{L''})>\olll(S_{L'',0})>
(1-\rho )\olll(F_{L''})>(1-3\gamma )2^{-M}.\end{equation}

By \eqref{**18b},
$X_{K+1,L''}(x)=(1-\rho ')\tCC_{\gamma }\cdot 1
=(1-\rho ')\tCC_{\gamma }\olll(F_{0})$ if $x\in S_{L'',0}.$

By \eqref{**18b} and \eqref{*21d}
if $x\in S_{L'',l}\sm S_{L'',l-1}$ then for $l=1,...,L''$,
\begin{equation}\label{**24b}
X_{K+1,L''}(x)=(1-\rho ')r_{0}\cdot \cdot \cdot r_{l}
\tCC_{\gamma }=(1-\rho ')\olll(F_{l})\tCC_{\gamma }.
\end{equation}
This
and \eqref{*30x}
imply that for $x\in S_{L'',L''}\sm S_{L'',L''-1}$
\begin{equation}\label{*27d}
X_{K+1,L''}(x)=(1-\rho ')\tCC_{\gamma }\olll(F_{L''})
<2^{-M}(1-\rho ')\tCC_{\gamma }<0.999\cdot 2^{-M+1}.
\end{equation}
Using \eqref{*h26} we have the following measure estimate:
$$\frac1{(1-\rho )r_{0}\cdot \cdot \cdot r_{l}}
\olll(S_{L'',0})>
\olll(S_{L'',l})>
(1-\rho )\frac1{r_{0}\cdot \cdot \cdot r_{l}}
\olll(S_{L'',0}),$$
which by \eqref{**h26} is equivalent to
\begin{equation}\label{**27c}\frac1{1-\rho }\olll(S_{L'',0})>
\olll(F_{l})\olll(S_{L'',l})>
(1-\rho )\olll(S_{L'',0}).
\end{equation}
Suppose for $l=0,...,M-1,$
$\ell '(l)$ is chosen so that
\begin{equation}\label{*28d}
X_{K+1,L''}(x)\geq 0.999\cdot 2^{-l}
\text{ when }x\in S_{L'',\ell '(l)},
\end{equation}
but $X_{K+1,L''}(x)< 0.999\cdot 2^{-l}$
for some $x\in S_{L'',\ell '(l)+1}$,
by \eqref{*27d} such an $\ell '(l)\leq L''$
exists.
By \eqref{**24b}
$$(1-\rho ')\olll(F_{\ell' (l)+1})\tCC_{\gamma }<
0.999\cdot 2^{-l},$$
and
$$(1-\rho ')\olll(F_{\ell' (l)})\tCC_{\gamma }\geq
0.999\cdot 2^{-l}$$
hold.
Therefore, using $\olll(F_{\ell '(l)+1})/\olll(F_{\ell '(l)})
>(1-2\gamma )$ we infer
\begin{equation}\label{*28c}
0.999\cdot 2^{-l}\leq
(1-\rho ')\olll(F_{\ell' (l)})\tCC_{\gamma }<\frac{0.999}{1-2\gamma }
2^{-l}.
\end{equation}
Set $S_{L'',\ell '(-1)}=\ess.$
By using \eqref{**24b} and the above definitions, estimates
for
$l=0,...,M-1,$
$x\in S_{L'',\ell '(l)}\sm S_{L'',\ell '(l-1)}$
we have
\begin{equation}\label{*31x}
0.999\cdot 2^{-l}\leq X_{K+1,L''}(x)<
0.999\cdot 2^{-(l-1)}.
\end{equation}
By \eqref{**27c}
$$\olll(S_{L'',\ell '(l)})<
\frac1{1-\rho }\olll(S_{L'',0})\frac1{\olll(F_{\ell' (l)})}<$$
(using \eqref{*27c} and \eqref{*28c})
$$\frac1{(1-\rho)^2 }\olll(F_{L''})\frac{(1-\rho ')\tCC_{\gamma }}
{0.999\cdot 2^{-l}}<
\frac{(1-\rho ')\tCC_{\gamma } 2^{l}}
{0.999\cdot (1-\rho)^2 }2^{-M},$$
on the other hand, by using \eqref{**27c}
$$\olll(S_{L'',\ell '(l)})>
\frac{(1-\rho )\olll(S_{L'',0})}
{\olll(F_{\ell '(l)})}>$$
(using \eqref{*27c} and \eqref{*28c} again)
$$
\frac{(1-\rho )^{2}\olll(F_{L''})(1-\rho ')\tCC_{\gamma }
(1-2\gamma )}
{0.999\cdot 2^{-l}}>$$
(using \eqref{*28j})
$$
 (1-\rho )^{2}(1-2\gamma )^{2}2^{-M}\cdot 2^{l}(1-\rho ')
\tCC_{\gamma }.
$$
Thus using \eqref{*108x}
for $l=0,...,M-1$
\begin{equation}\label{*33dfw}
\olll(S_{L'',\ell '(l)}\sm
S_{L'',\ell '(l-1)})>
\end{equation}
$$
 (1-\rho )^{2}(1-2\gamma )^{2}
2^{-M+l}
(1-\rho ')\tCC_{\gamma }-
\frac{(1-\rho ')\tCC_{\gamma }}
{0.999\cdot
(1-\rho)^2 }\cdot
2^{-M+l-1}>0.99\cdot 2^{-M+l-1}.
$$
By \eqref{*31x} if
$l=0,...,M-1,$
 $x\in S_{L'',\ell '(l)}\sm S_{L'',\ell '(l-1)}$
we have $X_{K+1,L''}(x)\geq 0.999\cdot 2^{-l}.$

By \eqref{**24b}, $X_{K+1,L''}$ takes different constant
values
on the set
$S_{L'',0}$ and on the sets $S_{L'',l}\sm S_{L'',l-1}$
for $l=1,...,L''.$

When $K>0$ we also know that $X_{K+1,L''}$
is pairwise independent  from $X_{h}$ for $h=1,...,K.$ Hence,
any function which is constant on the sets $S_{L'',0}$,
$S_{L'',l}\sm S_{L'',l-1}$, $l=1,...,L''$ is still
 independent from each $X_{h}$ for $h=1,...,K.$

Set $X_{K+1}'(x)=0.99\cdot 2^{-l}$ if $x\in S_{L'',\ell '(l)}\sm
S_{L'',\ell '(l-1)}$ for $l=0,...,M-1.$
Set $X_{K+1}'(x)=0$ if $x\not\in S_{L'',\ell '(M-1)}.$
Now $X_{K+1}'\leq X_{K+1,L''}$.
When  $K>0, X_{K+1}'$ is still independent
from each $X_{h}$, $h=1,...,K$ and 
it takes its values in $ \{ 0, 0.99, ..., 0.99\cdot 2^{-M+1}  \}$.
But it is $M$$-$$0.99$ super distributed.
By Lemma \ref{*super} we can choose an
$M$$-$$0.99$-distributed $X_{K+1}\leq X_{K+1}'$
 which is still independent from $X_{h}$ for each $h=1,...,K.$
This completes the part of our proof when we build
the $1-M$ family, that is for $K=0$ our argument ends here.

\subsubsection{The $K>0$ cases of our induction\\ Step $L=0$
of our leakage}\label{422b}

Next we assume that $K\geq 1$ and we can define
 $K-M$ families.

We use the definitions of the first paragraph of
Subsection \ref{422}. After the definition of $X_{K+1,0}$
we argue this way:

Choose a $K-M$ family on $\R$ with input constants
$\delta _{0}=\frac{\delta }{4(L'+1)}$,
$\OOO_{0}=\OOO$, $\GGG_{0} C_{\gamma }< \GGG$,
 $A_{0}=A$, $\capp_{0}$.
Then there exist a period $\tau_{0} $;
functions  $f_{h,0}:
\R\to [0,\oo)$,
pairwise independent  $M$$-$$0.99$-distributed ``random" variables
$X_{h,0}:\R\to\R$, for $h=1,...,K,$
a set $E_{\delta_{0} }$ periodic by $\tau_{0} ,$
 with
$\olll(E_{\delta_{0} })<\delta_{0} .$
Moreover, for all $x\not\in E_{\delta_{0} }$, there exist
$\omega_{0} (x)>\alpha_{0} (x)>A,$ $\tau_{0} (x)< \tau_{0} $
such that
$\omega_{0} ^{2}(x)<\tau_{0} ,$
$\frac{\omega_{0} (x)}{\alpha_{0} (x)}>\OOO \tau_{0} (x)
=\OOO_{0} \tau_{0} (x)$,
if $\alpha_{0} (x)\leq n<n+m\leq \omega_{0} (x)$,
and $\tau_{0} (x)|m$ then for all $h=1,...,K+1$
(for $h=1,...,K$ by the definition of the $K-M$ family,
for $h=K+1$ by the definition in the first
line of Subsection \ref{422}) there exists
$0\leq  f_{h,0}$ such that
$$
\frac{1}m\sum_{k=n}^{n+m-1}f_{h,0}(x+k^{2})>X_{h,0}(x).$$
For all $p\in\capp_{0}$, $(\tau_{0} (x),p)=1,$ $(\tau_{0} ,p)=1.$
For all $x\not\in E_{\delta_{0} }$ and all $h=1,...,K+1,$
$f_{h,0}(x+j+\tau_{0} (x))=f_{h,0}(x+j)$ whenever
$\alpha_{0} ^{2}(x)\leq j<j+\tau_{0} (x)\leq \omega_{0} ^{2}(x).$
Finally,
$$
\frac{1}{\tau_{0}} \int_{0}^{\tau_{0} }f_{h,0}=\ointt f_{h,0}<
C_{\gamma }\GGG_{0}\cdot
2^{-M+1},
$$
for $h=1,...,K.$

\subsubsection{Case $K>0$,
the setting after step $L-1$ of the leakage}
\label{423b}

We can repeat almost exactly the argument of Subsection
\ref{423}. We only need to add after the paragraph ending
with \eqref{*c17} that
 for $h=1,...,K$
\begin{equation}\label{*h11}
\frac{1}{\tau_{L-1}} \int_{0}^{\tau_{L-1} }f_{h,L-1}=\ointt f_{h,L-1}<
C_{\gamma }\GGG_{0}\cdot
2^{-M+1}.
\end{equation}
We emphasize that we do not expect that \eqref{*h11}
holds for $h=K+1$ and continue with the paragraphs
of Subsection \ref{423}
concerning the distribution of $X_{K+1,L-1}$.

\subsubsection{Case $K>0$, rearrangement with respect to $\tau '_{L-1}$,
choice of $\tau '_{L-1}$}
\label{424b}

This subsection is again almost completely identical to
Subsection \ref{424}. The only extra remark we need
after the first line of the last paragraph of
Subsection \ref{424} is the following:
We also have
\begin{equation}\label{*15}
\ointt f_{h,L-1}'\leq
\ointt f_{h,L-1}<C_{\gamma }\GGG_{0}\cdot 2^{-M+1},
\end{equation}
for $h=1,...,K.$

\subsubsection{Case $K>0$, choice of
$\kappa _{L}$,
$q_{L}$, $\FFF_{L}$, $\PPP_{L}$,
and
$X_{h,L}$ on $\R\sm\oLLL'(q_{L})$}\label{425b}

This subsection is identical to Subsection
\ref{425}.

\subsubsection{Case $K>0$,
putting $K-M$ families on $\oLLL'(q_{L})$}
\label{426b}

This is the subsection where we have a huge difference.
This is where we need to use the results from the previous
step of the induction on $K$.

The first four paragraphs until the definition
of $\ocap_L$ are identical to the ones in
Subsection \ref{426}.

Contrary to the $K=0$ case now we have to put
a $K-M$ family on $\oLLL'(q_{L})$.

For the choice of the $K-M$ family living on
$\oLLL'(q_{L})$ use Lemma \ref{lputres} with
$\ocap_{L},$ $\delta _{L}=\delta /4(L'+1)$,
$\OOO_{L}=\OOO\cdot q_{L}\tau _{L-1}'$,
$\Gamma_0$ and
$A.$

\begin{enumerate}
\item We obtain functions $\off_{h,L}$, $\oXX_{h,L}$
periodic by $\ottt_{L}q_{L}$ for $h=1,...,K$,
where $\ottt_{L}$ is a suitable natural number.
The functions $\oXX_{h,L}:\R\to\R$ are pairwise independent
and
conditionally $M$$-$$0.99$ distributed on $\oLLL'(q_{L}).$
There exists $E_{\delta _{L}}$ periodic by $\ottt_{L}q_{L}$.
For $h=1,...,K$ and $x\not\in \oXXX(q_{L})=
\cup_{j\in\Z}[j q_{L}, jq_{L}+\gggg 2^{\kkk_{L}})$ we have $\off_{h,L}(x)=0.$
\item
We have $\olll(E_{\delta _{L}})<\delta _{L}=\delta /4(L'+1)$.
For all $x\not\in E_{\delta _{L}}$,
there exist 
$\oooo_{L}(x)>\oaaa_{L}(x)>A$, $\ottt_{L}(x)< \ottt_{L}q_{L},$
$\oooo_{L}^{2}(x)<\ottt_{L}q_{L}$,
$$\frac{\oooo_{L}(x)}{\oaaa_{L}(x)}>\OOO_{L}
\ottt_{L}(x)=\OOO q_{L}\tau _{L-1}'\ottt_{L}(x).$$
Moreover,
if
$\oaaa_{L}(x)\leq n<n+m\leq\oooo_{L}(x)$
and $\ottt_{L}(x)|m$ then for $h=1,...,K$,
\begin{equation}\label{*18}
\frac{1}m\sum_{k=n}^{n+m-1}\off_{h,L}(x+k^{2})>
\oXX_{h,L}(x).
\end{equation}
\item
For all $p\in\ocap_{L},$
 $(\ottt_{L}(x),p)=1$,
$(\ottt_{L}q_{L},p)=1$.
\item
For all $x\in \oLLL'(q_{L})\sm E_{\delta _{L}},$
for all $h=1,...,K$,
$$\off_{h,L}(x+j+\ottt_{L}(x))=\off_{h,L}(x+j),$$
when $\oaaa_{L}^{2}(x)\leq j<j+\ottt_{L}(x)\leq \oooo^{2}_{L}(x)$.
\item
Finally, for all $h=1,...,K$
\begin{equation}\label{*25x}
\ointt\off_{h,L}<\GGG_{0}\cdot \gamma\cdot  2^{-M+1}.
\end{equation}
\end{enumerate}

We now define  the functions $X_{h,L}$
on $\oLLL'(q_{L})$ for $h=1,...,K$, by 
setting $X_{h,L}(x)=\oXX_{h,L}(x)$ if $x\in \oLLL'(q_{L}),$
$h=1,...,K$; also define
\begin{equation}\label{*19}
f_{h,L}=f'_{h,L-1}\cdot\chi _{\FFF_{L}}
+\off_{h,L} \text{ for }h=1,...,K.
\end{equation}
Where $f'_{h,L-1}$ is defined in Section \ref{424}
and $\FFF_{L}$ in Section \ref{425}.
It is important that by \eqref{*fi3}, $\oXXX(q_{L})$, which contains the support of
$\off_{h,L}$ is disjoint from 
$\FFF_{L}$ which contains the support of $f'_{h,L-1}\cdot\chi _{\FFF_{L}}$.

Recall
from \eqref{*15} that
$\ointt f'_{h,L-1}\leq \ointt f_{h,L-1}$
for $h=1,...,K.$

Using
\eqref{*25x} and
that $f'_{h,L-1}$ is periodic by $\tau '_{L-1}$
and $\FFF_{L}\sse \R\sm \oLLL'(q_{L})$ consists of blocks of length
$\tau '_{L-1}$
\begin{equation}\label{*fhlest1}
\ointt f_{h,L}\leq \olll(\FFF_{L})
\cdot \ointt f_{h,L-1}+
\GGG_{0}\left(\frac{\gamma }{
\olll(\oLLL'(q_{L}))}
\right)
2^{-M+1}\olll(\oLLL'(q_{L}))=(*).
\end{equation}
Since $\FFF_{L}\sse \R\sm \oLLL'(q_{L})$ by \eqref{*15}
we have
$$\left (\ointt f_{h,L-1}
\right)
\frac{\olll(\FFF_{L})}
{\olll(\R\sm \oLLL'(q_{L}))}<
\GGG_{0}\cdot C_{\gamma }2^{-M+1}.
$$
Hence,
we can continue our estimation by using \eqref{*080211}
$$
(*)<\left (\ointt f_{h,L-1}
\right) \frac{\olll(\FFF_{L})}
{\olll(\R\sm \oLLL'(q_{L}))}\cdot \olll(\R\sm
\oLLL'(q_{L}))+
\GGG_{0} C_{\gamma }2^{-M+1}\olll(\oLLL'(q_{L}))<$$
\begin{equation}\label{*fhlest2}
 \GGG_{0} \cdot C_{\gamma }\cdot 2^{-M+1}.
 \end{equation}

After these observations we can return to Section \ref{426},
to the definition of $\capp_{L}$ and read everything
until the paragraph marked by a $\dagger$.

When $K>0$ we need to add the following estimate to the case
when $x\in \R\sm (\oLLL'(q_{L})\cup E^{L})\sse \tFFF_{L}\sse \FFF_{L}$.

If $\alpha _{L}(x)\leq n<n+m\leq\omega _{L}(x)$
and
$\tau _{L}(x)|m$ then by \eqref{*f12} and \eqref{*18b}
 for $h=1,...,K$
\begin{equation}\label{*12x}
\frac1{m}\sum_{k=n}^{n+m-1}f'_{h,L-1}(x+k^{2})
=\frac1{m}\sum_{k=n}^{n+m-1}f_{h,L}(x+k^{2})
>X_{h,L}(x)=X_{h,L-1}'(x).
\end{equation}

Next suppose $x\in \oLLL'(q_{L})\sm E^{L}\sse\oLLL'(q_{L})\sm
E_{\delta _L}$. For $h=1,...,K$
the estimates which we have for the $K-M$ family
put on $\oLLL'(q_{L})$ can be applied.
In other words, for these $x$,
set
$\omega _{L}(x)\defeq \oooo_{L}(x)>
\alpha _{L}(x)\defeq \oaaa_{L}(x),$
$\tau _{L}(x)\defeq \ottt_{L}(x)q_{L}\tau _{L-1}' <\tau _{L}=q_{L}\ottt_{L}\tau'_{L-1}$.
Then
$\omega_{L}^{2}(x)<\tau _{L} $, and
$$\frac{\omega _{L}(x)}{\alpha _{L}(x)}>\OOO_{L}\ottt_{L}(x)=
\OOO\tau _{L}(x)
.$$
Furthermore, if $\alpha _{L}(x)\leq n<n+m\leq\omega _{L}(x)$
and $\tau _{L}(x)|m$, then
$\ottt_{L}(x)|m$ and
by \eqref{*18} and \eqref{*19} 
we have for $h=1,...,K$
$$\frac1{m}\sum_{k=n}^{n+m-1}f_{h,L}(x+k^{2})
>X_{h,L}(x).$$
For all $p\in \capp_{L}\sse \ocap_{L}$ we have
$(\ottt_{L}(x)q_{L}\tau '_{L-1},p)=(\tau _{L}(x),p)=1$ and
for
$h=1,...,K$, if
$\alpha _{L}^{2}(x)\leq j<j+\tau _{L}(x)\leq \omega _{L}^{2}(x)$
then $\alpha _{L}^{2}(x)=\oaaa_{L}^{2}(x)\leq j<j+\ottt_{L}(x)<
...<j+q_{L}\tau '_{L-1}\ottt_{L}(x)\leq
 \omega _{L}^{2}(x)=\oooo_{L}^{2}(x)$
and hence
$
f_{h,L}(x+j+\tau _{L}(x))=
\off_{h,L}(x+j+\tau_{L}(x))=
\off_{h,L}(x+j+q_{L}\tau '_{L-1}\ottt_{L}(x))=
\off_{h,L}(x+j+(q_{L}\tau _{L-1}'-1)\ottt_{L}(x))=...=
\off_{h,L}(x+j)=
f_{h,L}(x+j).$

\subsubsection{Case $K>0$, properties of $f_{K+1,L}$}
\label{427b}

This section is again identical to Section \ref{427}.

\subsubsection{Case $K>0$, finishing the leakage}\label{428b}

We start to argue as in Subsection \ref{428}.
We need to insert just before the sentence
containing \eqref{*24b} the remark:
Moreover, if $\alpha (x)\leq n<n+m\leq \omega (x)$
and
$\tau (x)|m$, then for all $h=1,...,K$
letting $f_{h}= f_{h,L''}$
(see also \eqref{*c17} which is used with $L=L''+1$)
\begin{equation}\label{*h44}
\frac{1}m\sum_{k=n}^{n+m-1}f_{h}(x+k^{2})>X_{h}(x).
\end{equation}

Before \eqref{*14bx} we need to add the comment
that by (\ref{*fhlest1}-\ref{*fhlest2}) for $h=1,...,K$
\begin{equation}\label{*48yy}
\ointt f_{h}<
C_{\gamma }\GGG_{0}2^{-M+1}<
\GGG\cdot
2^{-M+1}.
\end{equation}

The rest of the argument is identical to Subsection
\ref{428}
and this way we can complete our induction.

\end{proof}

\section{Proof of the Main Result}
\label{secmain}

Lemma \ref{lksystem} yields the next theorem which,
as we will see, easily implies Theorem \ref{mainth}.

\begin{theorem}\label{*c1c}
Given $\delta >0$, $M$ and $K$ there exist
$\tau_0 \in \N,$ $\oEE_{\delta }\sse [0,1),$
a measurable transformation $T:[0,1)\to[0,1)$,
$T(x)=x+\frac{1}{\tau_0}$ modulo $1$,
$f:[0,1)\to [0,+\oo)$,
$\oXX_{h},$ $h=1,...,K$ which are
pairwise independent  $M$$-$$0.99$-distributed
random variables defined on $[0,1)$ equipped with
the Lebesgue measure, $\lambda $, such that
$\lambda (\oEE_{\delta })<\delta $, for all
$x\in [0,1)\sm \oEE_{\delta }$ there exists
$N_{x}$ satisfying
$$
\frac{1}{N_{x}}\sum_{k=1}^{N_{x}}
f(T^{k^{2}}(x))>\sum_{h=1}^{K}\oXX_{h}(x),
$$
and $\int_{[0,1)}fd\lambda <K\cdot 2^{-M+2}.$
\end{theorem}

\begin{proof}
Use Lemma \ref{lksystem} with $\delta $, $\OOO=1000,$
$\GGG=1.1$, $A=1,$ $\capp=\ess$ to obtain a $K-M$ family
with $E_{\delta }$, $f_{h}$ and $X_{h}$ periodic by $\tau =\tau_0 .$
Set $\oEE_{\delta }=\frac{1}{\tau_0 } E_{\delta }\cap [0,1)$
and
for $x\in [0,1)$ set $\off_{h}(x)= 1.01\cdot f_{h}(\tau_0 \cdot x),$
$\oXX_{h}(x)=X_{h}(\tau_0 \cdot x)$.

Assume
$x\in [0,1)\sm \oEE_{\delta }$.
Since $\OOO \alpha (\tau_0\cdot  x)
\tau (\tau_0 \cdot x)<\omega (\tau_0 \cdot x)$ we have
$\alpha (\tau_0 \cdot x)=n<n+(\OOO-1)\alpha (\tau_0 \cdot x)\tau (\tau_0 \cdot x)
<\omega (\tau_0\cdot  x)$ and \eqref{*12dfw} used
with $n=\aaa(\tttt_{0}\cdot x)$ and $m=(\OOO-1)\aaa(\tttt_{0}\cdot x)\tttt(
\tttt_{0}\cdot x)$ implies
$$\frac{1}{(\OOO-1)\alpha (\tau_0 \cdot x) \tau (\tau_0 \cdot x)}
\sum_{k=\alpha (\tau_0 \cdot x)}^{
\alpha (\tau_0 \cdot x)+(\OOO-1)\alpha (\tau_0
\cdot x) \tau (\tau_0 \cdot x)-1}
f_{h}(\tau_0 \cdot (T^{k^{2}}x))
=$$
$$\frac{1}{(\OOO-1)\alpha (\tau_0 \cdot x)\tau (\tau_0 \cdot x)}
\sum_{k=\alpha (\tau_0 \cdot x)}^{
\alpha (\tau_0 \cdot x)+(\OOO-1)\alpha (\tau_0 \cdot x) \tau (\tau_0 \cdot x)-1}
f_{h}(\tau_0 \cdot x+k^{2})>X_{h}(\tau_0 \cdot x).
$$
Since $f_{h}\geq 0$, if we let $N_{x}=
\alpha (\tau_0 \cdot x)+(\OOO-1)\alpha (\tau_0
\cdot x) \tau (\tau_0 \cdot x)-1$,
then  since $\OOO=1000$, $N_{x}/(\OOO-1)\aaa(\tttt_{0}\cdot x)\tttt(
\tttt_{0}\cdot x)<1.01$, for all $h=1,...,K$
\begin{equation}\label{*fhbar}
\frac{1}{N_{x}}\sum_{k=1}^{N_{x}}
\off_{h}(T^{k^{2}}x)
=
\frac{1.01}{N_{x}}\sum_{k=1}^{N_{x}}
f_{h}(\tau_0 \cdot (T^{k^{2}}x))
\geq \oXX_{h}(x).
\end{equation}
Let $f(x)$ be the restriction
of $\sum_{h=1}^{K} \off_{h}(x)$
onto $[0,1).$
Therefore,
using \eqref{*17iest} with $\gggg'=1$
from Lemma \ref{lksystem} we obtain
 $$\int_{0}^{1}f(x)d\lambda (x)
= 1.01\sum_{h=1}^{K}\ointt f_{h}<1.01\cdot \GGG\cdot K \cdot
2^{-M+1}<
K\cdot 2^{-M+2}.$$
For all $x\in [0,1)\sm \oEE_{\delta }$ by \eqref{*fhbar} there exists
$N_{x}$ such that
$$\frac{1}{N_{x}}\sum_{k=1}^{N_{x}}
f(T^{k^{2}}x)>\sum_{h=1}^{K}\oXX_{h}(x).$$
\end{proof}

\medskip

Now we can complete the proof of Theorem \ref{mainth}.

\begin{proof}
For each $p\in\N$ set $M_{p}=4^{p}.$
On the probability space
$([0,1),\lambda )$
consider $M_{p}-0.99$-distributed
random variables $\oXX_{h}$ for $h=1,...,K$ for a sufficiently
large $K$. Assume that $u$ denotes the mean of these variables.
An easy calculation shows that
$$u=\int_{[0,1]} \oXX_{h}(x)d\lambda (x)=
\sum_{l=0}^{M_{p}-1}0.99^{2}\cdot 2^{-l}\cdot 2^{-M_{p}+l-1}>
0.9\cdot M_{p}\cdot 2^{-M_{p}-1}.$$
By the weak law of large numbers
$$\lambda \left  \{ x:\left  |
\frac{1}K\sum_{h=1}^{K}\oXX_{h}(x)-u
\right|\geq \frac{u}2\right  \}\to 0.$$
Fix $K$ so large that
$$\lambda \left  \{ x:\frac{1}K
\sum_{h=1}^{K}
\oXX_{h}(x)\geq \frac{u}2
\right  \}>1-\frac{1}p,$$
and let
$$U_{p}'=\left  \{ x:
\frac{1}K\sum_{h=1}^{K}\oXX_{h}(x)>\frac{0.9}2
\cdot M_{p}\cdot 2^{-M_{p}-1}\right  \}.$$
We have $\lambda (U_{p}')>1-\frac{1}p.$
By Theorem \ref{*c1c} used with $\delta =\frac{1}p,$
$M_{p}$ and $K$ there exist $\tau_0 \in \N,$
$\oEE_{1/p}\sse [0,1)$ and a periodic transformation
$T:[0,1)\to [0,1),$ $T(x)=x+\frac{1}{\tau_0 }$ modulo $1,$
$f:[0,1)\to [0,+\oo)$,
$\oXX_{h}$ pairwise independent
$M_{p}-0.99$-distributed random variables defined on
$[0,1)$ such that $\lambda (\oEE_{1/p})<\frac{1}p$
and for all $x\in [0,1)\sm \oEE_{1/p}$
there exists $N_{x}$ such that
$$\frac{1}{N_{x}}\sum_{k=1}^{N_{x}}
f(T^{k^{2}}x)>\sum_{h=1}^{K}\oXX_{h}(x)$$
and
$\int_{[0,1)}fd\lambda <K\cdot 2^{-M_{p}+2}.$
Put $U_{p}=U_{p}'\sm
\oEE_{1/p}.$
Then $\lambda (U_{p})>1-\frac{2}p$ and for
$x\in U_{p}$ there exists $N_{x}$ such that
$$\frac{1}{N_{x}}
\sum_{k=1}^{N_{x}}f(T^{k^{2}}x)>\sum_{h=1}^{K}
\oXX_{h}(x)>K\cdot \frac{0.9}2\cdot M_{p}\cdot 2^{-M_p-1}.$$
Thus letting $t_{p}=K\cdot \frac{0.9}2\cdot M_{p}\cdot 2^{-M_p-1}$,
and $$\tUU_{p}=\left  \{ x:\sup_{N}
\frac{1}N
\sum_{k=1}^{N}f(T^{k^{2}}x)>t_{p}\right  \}$$
we have $U_{p}\sse \tUU_{p}$ and hence
$\lambda (\tUU_{p})>1-\frac{2}p.$
On the other hand
$$\frac{\int f d\lambda }{t_{p}}=\frac{\int |f|d\lambda }{t_{p}}<
\frac{K\cdot 2^{-M_{p}+2}}
{K\cdot \frac{0.9}2\cdot M_{p}\cdot 2^{-M_p-1}}
< \frac{32}{M_{p}}.$$
Hence, $\lambda (\tUU_{p})\to 1$ and $\int |f|d\lambda/t_{p} \to 0$
as $p\to \oo.$
Therefore there is no $C$ for which \eqref{weaklolo}
holds with $\mu=\lambda$.
This implies that the sequence $n_{k}=k^{2}$ is
$L^{1}$-universally bad.
\end{proof}

 {\bf ACKNOWLEGMENT} The authors thank an anonymous referee for
pointing out the use of the P\'olya-Vinogradov inequality in Lemma 3.
It shortened our orginal argument.
Other comments of the referees also helped to improve the paper.

\end{document}